\documentclass[oneside]{amsart}
\usepackage{amsmath,ifthen,amscd,amssymb,
graphicx,subfigure}
\newcommand{\NSgeomfinite}{MR96c:20066}
\newcommand{\Enonhopf}{MR2022477}
\newcommand{\Ethesis}{Ethesis}
\newcommand{\Wise}{MR96m:20058}

\newcommand{\Cannon}{MR88a:20049}

\newcommand{\BB}{BB}
\newcommand{\BH}{MR1744486}

\newcommand{\HU}{MR83j:68002}
\newcommand{\Epstein}{MR93i:20036}
\newcommand{\LS}{MR58:28182}
\newcommand{\Gersten}{MR94j:20043}

\graphicspath{{.}{./eps/}}

\newtheorem{thm}{Theorem}[section]
\newtheorem{lem}[thm]{Lemma}

\newtheorem{cor}[thm]{Corollary} 
\newtheorem{conj}[thm]{Conjecture}
\theoremstyle{definition}
\newtheorem{defn}[thm]{Definition}

\newtheorem{eg}{Example}

\newcommand{\blackboard}[1]{\ensuremath{\mathbb{#1}}}

\newcommand{\smallcaps}[1]{\textrm{\textsc{#1}}}

\newcommand{\N}{\blackboard{N}}
\newcommand{\Z}{\blackboard{Z}}

\newcommand{\cat}{\smallcaps{CAT}}

\newcommand{\iif}{isoperimetric function}
\newcommand{\ac}{almost convex}
\newcommand{\fftp}{falsification by fellow traveler property}
\newcommand{\cg}{Cayley graph}
\newcommand{\hnn}{HNN extension}

\newcommand{\gset}{generating set}

\newcommand{\bbx}{$(G_{1,1},X)$}
\newcommand{\wx}{$(G_W,Y)$}

\newcommand{\fsa}{finite state automaton}
\newcommand{\mt}{\mathtt}
\newcommand{\ra}{\rightarrow}
\begin{document}

\title[Patterns theory and geodesic automatic structure]
  {Patterns theory and geodesic automatic structure for a class of groups}

\author[M.~Elder]{Murray J. Elder}
\address{(As at Dec 2006) Dept of Mathematical Sciences\\
	Stevens Institute of Technolgy\\
	Hoboken NJ 07030 USA}
\email{melder@stevens.edu}

\date{Published in  International Journal of Algebra and Computation Vol 13, No 2 (2003) 203--230}

\begin{abstract}
We introduce a theory of {\em patterns} in order to study geodesics in
a certain class of group presentations. Using patterns we show that
there does not exist a geodesic automatic structure for certain group
presentations, and that certain group presentations are almost convex.
\end{abstract}

\maketitle

\section{Introduction}
In this article we examine the possibility of automatic structures in
two intriguing examples.  Wise's example is non-Hopfian and \cat(0)
\cite{\Wise}, so proving this group is not automatic would answer the
open question: Is every \cat(0) group automatic?  Conversely proving
this group is automatic would answer another open question: Is every
automatic group Hopfian (or residually finite)?  Brady and Bridson's
example is not biautomatic, not \cat(0), and has a quadratic \iif\
\cite{\BB}. Finding an automatic structure for this group would answer
a third open question of fundamental interest to automatic group
theorists: Does every automatic group admit some biautomatic
structure?

After some preliminary definitions we examine in some detail Brady and
 Bridson's example with a chosen \gset, which we denote by the pair
 \bbx.
We define the notion of a {\em pattern} in the \cg\ and this idea
 enables us to prove that for the chosen \gset\ the full language of
 geodesics is not regular (Theorem \ref{notreg}), and further that for
 the chosen \gset\ the group can have no geodesic automatic structure
 (Corollary \ref{nogeodauto}).  We characterize the set of all
 patterns for \bbx\ (Theorem \ref{patternschar}), and use this to
 prove that the pair is \ac\ (Theorem \ref{acthm}).  We remark that
 all of the above results can be replicated for Wise's example.

While these results go nowhere close to resolving the question of
 automaticity for these two examples, they provide much insight into
 the geodesic structure of this class of groups.

The author wishes to thank his advisor Walter Neumann for his
encouragement and ideas throughout this project.

\section{Preliminaries: Automatic Groups}
Let $(G,X)$ be a group with finite \gset\ $X$, assume that $X$ is
inverse closed, and let $X^*$ denote the set of all possible words in
the letters of $X$, including the empty word.  Let $\Gamma(G,X)$ be
the \cg\ for the pair $(G,X)$, metrized by giving each edge unit
length and endowing it with the path metric. That is, the distance
between any two points is the infimum of the lengths of paths between
these points in the graph.  A word in $X^*$ describes a path in the
\cg\ and vice versa.  Paths can be parameterized by non-negative $t\in
\mathbb {R}$ by defining $w(t)$ as the point distance $t$ along the
path $w$ from the identity if $t<|w|$ and $w(t)=\overline w$ if
$t\geq|w|$, where $\overline w$ is the endpoint of $w$ and $|w|$ is
the length of the path $w$ which is equal to the number of letters in
the word $w$.  If two words $w,u$ evaluate to the same group element
in $G$ then we write $w=_Gu$.

Paths $w$ and $u$ are said to {\em $k$-fellow travel} if $d(w(t),
u(t)) \leq k $ for each $ t\in \mathbb {R}$ with $t\geq 0$.  The two
paths are {\em asynchronous} $k$-fellow travelers if there is a
non-decreasing proper continuous function $\phi : [0,\infty)
\rightarrow [0,\infty)$ such that $d(w(t), u(\phi (t))) \leq k $.
This means that any point on $w$ is within $k$ of some point on $u$
and vice versa.  We imagine the two paths traveling at different
speeds (but not backtracking) to keep within $k$ of each other.

\begin{defn}[Fellow traveler property]
A language $L\subset X^*$ enjoys the {\em (asynchronous) fellow
traveler property} if there is a constant $k$ such that for each
$w,u\in L$ with $d(\overline w,\overline u) \leq 1$ in $\Gamma(G,X)$,
$w$ and $u$ (asynchronously) $k$-fellow travel.
\end{defn}

A language $L\subset X^*$ on any finite alphabet $X$ is said to be
{\em regular} if it is the set of words accepted by some \fsa. See
\cite{\Epstein} for details.

A simple yet powerful result about regular languages is the following.
\begin{lem}[Pumping Lemma]
Let $M$ be a \fsa\ on an alphabet $X$, having $n$ states.  If $w\in
L(M)$ is a word of length greater than $n$ then we can write $w=uzv$
with $|z|>0$ and $uz^iv \in L(M) \; \forall i\geq 0$.
\end{lem}

\noindent
\begin{proof}
If $|w| >n$ then as $M$ reads $w$ it must pass through the same state
more than once.  Let $u$ be the prefix of $w$ until it reaches this
state. Let $z$ be the next part of $w$ until it gets back to the
repeated state. Then $v$ is the remaining part of $w$.
%
Now since $w=uzv$ ends in an accept state so does $uz^iv $ for all
$i\geq 0$.
\end{proof}

The set of regular languages is known as a {\em formal language}.  In
general the set of words on an alphabet that are accepted by some kind
of computing machine is a formal language.  Formal language classes
can be arranged in a hierarchy known as the ``Chomsky hierarchy'', in
order of increasing complexity. The complexity of a language can be
seen in terms of the complexity of the machine which accepts it. Since
\fsa\ are the simplest type of computing machines the class of regular
languages is at the start of the Chomsky hierarchy.  
 More details can be found in \cite{\HU}.

\begin{defn}[Automatic]
Let $(G,X)$ be a group with finite \gset\ $X$, and let $L\subset X^*$
be some set of words on $X^*$ such that
\begin{enumerate}
\item $L$ is regular
\item $L$ has the (asynchronous) fellow traveler property
\item $L$ surjects to $G$.
\end{enumerate}
Then $L$ is an {\em (asynchronous) automatic structure} for  $(G,X)$.
\end{defn}
If $G$ has an automatic structure with respect to one \gset\ then it
has an automatic structure with respect to every \gset\
\cite{\Epstein} so we say that $G$ is {\em(asynchronously) automatic}
if it has an (asynchronous) automatic structure for some \gset.

\begin{defn}[Biautomatic]
$L\subset X^*$ is a {\em biautomatic structure} for  $(G,X)$
if $L$ and $L^{-1}$ are both  automatic structures for   $(G,X)$.
\end{defn}

\noindent If $L$ is an automatic structure then the only thing to
check for biautomaticity is that $L^{-1}$ has the fellow traveler
property.

If $G$ is automatic then it has at most a quadratic \iif\
\cite{\Epstein}.

If $L\subset X^*$ is a language surjecting to $G$ and has the fellow
traveler property then we say $L$ is a {\em combing} for $G$. If
$L^{-1}$ also has the fellow traveler property then we say $L$ is a
{\em bicombing}.

A group is {\em CAT(0)} if it acts properly discontinuously
co-compactly by isometries on some \cat(0) metric space (See
\cite{\BH}). If $G$ is \cat(0) then 
space on which is acts admits a bicombing by geodesics. This does not
prove that \cat(0) groups are biautomatic or even automatic since the
geodesic bicombing may not project onto the 1-skeleton and also be a
regular language. Therefore whether \cat(0) groups are (bi)automatic
is presently an open problem.  
\begin{eg}
Let \bbx\ be the group with presentation
$$ \langle a,b,c,d,s,t | c=ab=ba, d=ab^{-1}, s^{-1}as=c, t^{-1}at=d \rangle.$$
$G_{1,1}$ has no biautomatic structure, has a quadratic \iif\
\cite{\BB}, and is not \cat(0) \cite{\Gersten}.  It is asynchronously
automatic \cite{\Ethesis} and the proof of this is omitted.
\end{eg}

\begin{eg}
Let \wx\ be the group with presentation
$$\langle a,b,c,d,s,t | c=ab=ba, d=c^2, s^{-1}as=d, t^{-1}bt=d
\rangle.$$
$G_W$ is non-Hopfian, \cat(0) \cite{\Wise} and asynchronously
automatic \cite{\Ethesis}.
\end{eg}

\section{Preliminaries: \hnn s}
\begin{defn}[Multiple \hnn] \label{mult}
Let $(A,Z)$ be a group with finite \gset\ $Z$ and relations $R$, let
 $U_1,\ldots, U_n, V_1, \ldots, V_n$ be subgroups of $A$ and let
 $\phi_i:U_i \ra V_i$ be an isomorphism for each $i$.  The group
 $(G,X)$ with presentation
$$\langle Z,s_1,\ldots, s_n | R, s_i^{-1}u_is_i=\phi_i(u_i) \; \forall
u_i\in U_i, \; \forall i \rangle$$
is a {\em multiple \hnn} of $(A,Z)$. The generators $s_i$ are called
{\em stable letters}, and the pairs of $U_i,V_i$ are called {\em
associated subgroups}.
\end{defn}
If each $U_i$ is finitely generated by $\{ u_{i_j} \}$ and $\phi_i
(u_{i_j})=v_{i_j}$ then $V_i$ is finitely generated by $\{v_{i_j}\}$.
Thus $(G,X)$ has the finite presentation
$$\langle Z,s_1,\ldots, s_n | R, s_i^{-1}u_{i_j}s_i=v_{i_j} \; \forall
i, \; \forall j \rangle.$$

\begin{thm}[Britton's Lemma]
Let $(G,X)$ be a multiple \hnn\ with the presentation in Definition
\ref{mult} above. If $w\in X^*$ is freely reduced and $w=_G1$ then $w$
contains a sub-word of the form $s_i^{-1}u_{i_j}s_i$ or
$s_iv_{i_j}s_i^{-1}$ for some non-trivial $u_{i_j}\in U_i$ or $
v_{i_j}\in V_i$.
\end{thm}

\noindent
A proof can be found in \cite{\LS}.  We call the sub-word
$s_i^{-1}u_{i_j}s_i$ or $s_iv_{i_j}s_i^{-1}$ a {\em pinch}.
\begin{defn}
A word that does not admit any pinches is called {\em stable letter
reduced}.  If two words have the same sequence of stable letters then
we say they have {\em parallel stable letter structure}.
\end{defn}

\begin{eg}
If $n=1$ and $A=\Z$ then $(G,X)$ is a Baumslag-Solitar group. If
$A=\Z^r$ for some $r>1$ and $n>1$ then $(G,X)$ is sometimes called a
{\em generalized Baumslag-Solitar group}.
\end{eg}

\begin{eg}
$G_{1,1}$ is a double \hnn\ of $\Z^2\cong \langle a,b,c,d | c=ab=ba,
d=ab^{-1} \rangle$ associating cyclic subgroups $\langle a
\rangle,\langle c \rangle$ and $\langle a \rangle,\langle d \rangle$.
\end{eg}

\begin{eg}
$G_W$ is a double \hnn\ of $\Z^2\cong \langle a,b,c,d | c=ab=ba, d=c^2
\rangle$ associating cyclic subgroups $\langle a \rangle,\langle d
\rangle$ and $\langle b \rangle,\langle d \rangle$.
\end{eg}

\begin{figure}[ht!]
  \begin{center}
      \includegraphics[width=12cm]{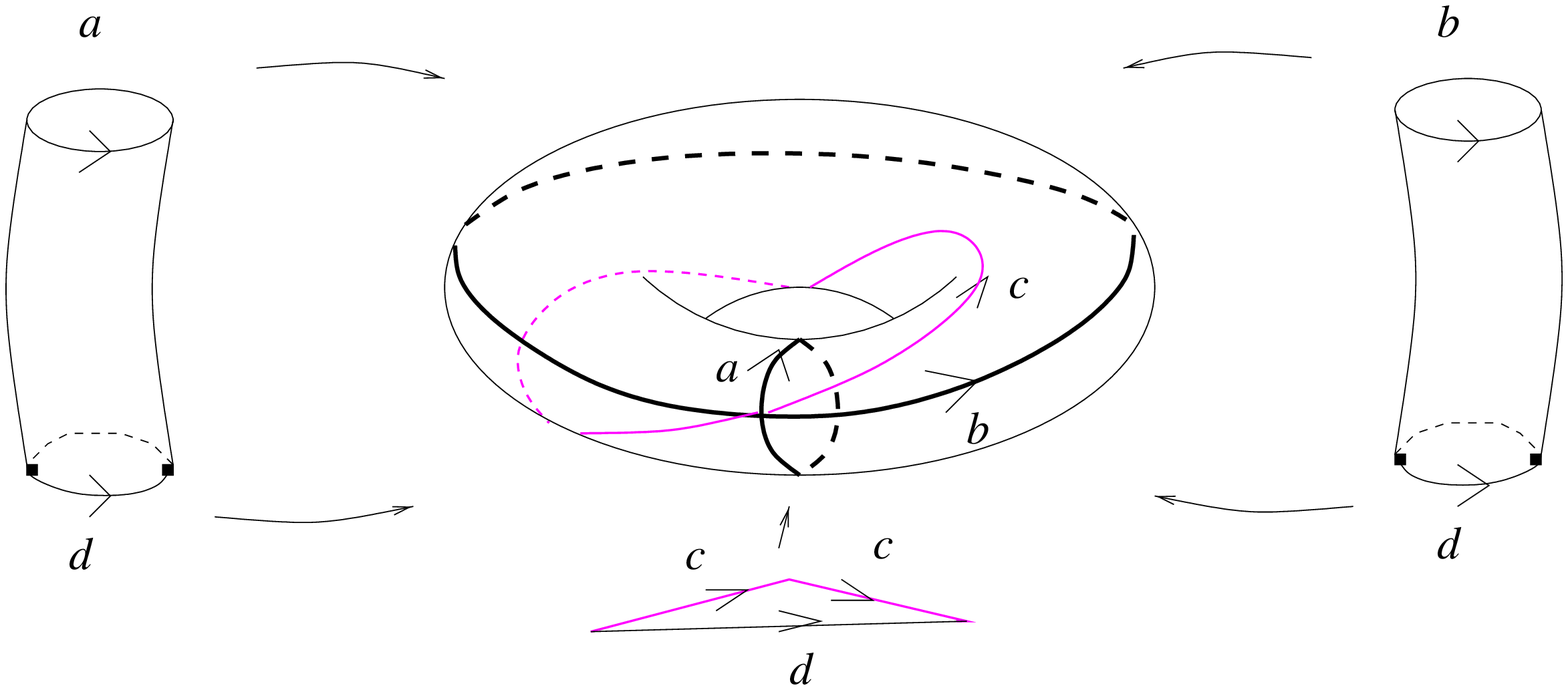}
  \end{center}
  \caption{Presentation 2-complex for \wx}
  \label{fig:pres2complexWX}
\end{figure}

The \cg\ of a multiple \hnn\ can be viewed in the following way.  For
each presentation of a group we can construct the ``presentation
2-complex'', where the 1-skeleton of its universal cover is the \cg,
and the 2-skeleton is sometimes called the ``Cayley complex'' or the
``filled \cg''.  If $(G,X)$ is a multiple \hnn\ of $(A,Z)$, the
presentation 2-complex is formed by taking the presentation 2-complex
for $(A,Z)$ and attaching $n$ annuli by appropriate edge gluings.
\begin{figure}[ht!]
  \begin{center}
      \includegraphics[width=10cm]{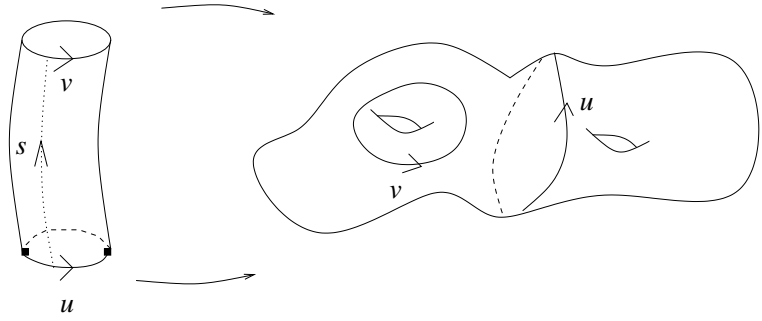}
  \end{center}
  \caption{Presentation 2-complex for a multiple \hnn}
  \label{pres2complex}
\end{figure}
In the universal cover we see copies of the ``base space'' or the
2-complex for $(A,Z)$ glued together by copies of the annuli, which
attach themselves along associated subgroups of $(A,Z)$.  For example,
the presentation 2-complex for \wx\ can be seen as two annuli (or
cylinders) with boundaries labeled $a,d$ and $b,d$ respectively,
attached to a torus with a triangle $c^2d^{-1}$ glued to it, by
identifying the loops $a,b$ and $d$, as in Figure
\ref{fig:pres2complexWX}.
The universal cover of the torus with a triangle attached is a plane
 with ``bumps'' on both sides, and the universal cover of each
 cylinder gives a bi-infinite ``strip''.  So we can view any such
 multiple \hnn\ as being made up of base-group {\em planes} glued
 together along associated subgroups by stable letter {\em strips}.

\begin{defn}[Strip equidistant]
Let $(G,X)$ be a multiple \hnn\ with the presentation in Definition
\ref{mult} above.  If $|u_i|=|\phi(u_i)|$ for all $i$ then we say
$(G,X)$ has a {\em strip equidistant} presentation.
\end{defn}
Note that if $(G,X)$ has a strip equidistant presentation then if a
word $w\in X^*$ admits a pinch then it can be shortened by $2$, so
geodesics are stable letter reduced.

\begin{defn}[Geodesic subgroup]
Let $(G,X)$ be any group with \gset\ $X$.  A subgroup $A$ of $G$ with
\gset\ $Y\subseteq X^*$ is {\em geodesic} in $(G,X)$ if each element
of $A$ has a geodesic word $w\in Y^*$ evaluating to it.
\end{defn}

\begin{defn}[Totally geodesic]
Let $(G,X)$ be any group with \gset\ $X$.  A subgroup $A$ of $G$ with
\gset\ $Y\subseteq X^*$ is {\em totally geodesic} in $(G,X)$ if every
geodesic word $w\in X^*$ evaluating to an element of $A$ is an element
of $Y^*$.  A subspace of the \cg\ which is a copy of $\Gamma(A,Y)$
inside $\Gamma(G,X)$ is called {\em totally geodesic} if $(A,Y)$ is.
\end{defn}

\noindent
It is clear that totally geodesic subgroups are geodesic.

\begin{eg}
\label{tot1}
\bbx\ has a strip equidistant presentation, and the subgroup $\Z^2
 \cong \langle a,b,c,d | ab=ba,c=ab,d=ab^{-1}\rangle$ is totally
 geodesic.  The associated subgroup $\langle a \rangle$ of $\Z^2$ is
 not totally geodesic in $(\Z^2, \{a,b,c,d\})$, since
 $c^id^i=_{Z^2}a^{2i}$. The subgroups $\langle a \rangle,\langle c
 \rangle,\langle d \rangle$ are geodesic in $(\Z^2, \{a,b,c,d\})$.

\end{eg}

\begin{eg}
\label{tot2}
\wx\ has a strip equidistant presentation, the subgroup $\Z^2 \cong
 \langle a,b,c,d | ab=ba,c=ab,d=c^{2}\rangle$ is totally geodesic in
 it, and all associated cyclic subgroups of $\Z^2$ are totally
 geodesic in $(\Z^2, \{a,b,c,d\})$.
\end{eg}

\section{Geodesics and sequences}
Throughout this section $(G,X)$ will be a multiple \hnn\ of $(\Z^n,Y)$
with a strip equidistant presentation, where $Y$ is some finite \gset\
for $\Z^n$.  Let a {\em $\Z^n$-plane} refer to the copies of
$\Gamma(\Z^n,Y)$ in $\Gamma(G,X)$.

\begin{lem}
$\Z^n$-planes are totally geodesic in $\Gamma(G,X)$.
\label{lem:totallygeod}
\end{lem}

\noindent
\begin{proof}
Suppose $w$ is a geodesic from $p$ to $q$ in a $\Z^n$-plane, and $w$
does not lie in the plane. Then $w$ must contain a stable letter. Let
$u$ be a path in the $\Z^n$-plane from $p$ to $q$. Then $wu^{-1}=_G 1$
so by Britton's Lemma it contains a pinch. Since $u$ has no stable
letters then $w$ must contain the pinch. But since the presentation is
strip equidistant performing a pinch would shorten $w$ by 2, which is
a contradiction.
\end{proof}

\begin{cor}
A geodesic path in  $\Gamma(G,X)$ visits no $\Z^n$-plane twice.
\end{cor}

\noindent
\begin{proof}
If a path leaves some $\Z^n$-plane at a point $p$ then re-enters at a
 point $q$, then it cannot be geodesic by the Lemma.
\end{proof}

\begin{lem}
Let $g\in G$. Every stable letter reduced path from $1$ to $g$ crosses
the same succession of planes and strips.
\label{lem:parallel}
\end{lem}

\noindent
\begin{proof}
We proceed by induction on the number of stable letters in a stable
letter reduced word.  If $w$ is a stable letter reduced and has no
stable letters, and if $u$ is a stable letter reduced word for $w$,
then $wu^{-1}=_G 1$ so if $wu^{-1}$ contains any stable letters it
must contain a pinch. Thus $u$ has no stable letters so both paths
stay in the base $\Z^n$-plane.  Now assume the hypothesis for stable
letter reduced words having at most $k-1$ stable letters. Suppose $w$
is stable letter reduced and has $k$ stable letters, and $u$ is a
stable letter reduced word for $w$. Write $w=w_0 r_1 \ldots w_{n-1}
r_n w_n$ and $u=u_0 t_1 \ldots u_{m-1} t_m u_m$.  Now $wu^{-1}=_G 1$
and contains stable letters so there is a pinch, so choose $i,j$
maximal such that $r_i(w_i \ldots r_n w_n u_m^{-1} t_m^{-1} \ldots
u_j^{-1})t_j^{-1}$ is a pinch.  Then $(w_i \ldots r_n w_n u_m^{-1}
t_m^{-1} \ldots u_j^{-1})$ is a subgroup element so if it contains a
stable letter it must contain a pinch, so $i,j$ were not maximal.

Thus $ r_n w_n u_m^{-1} t_m^{-1}$ is a pinch. Let $y= r_n w_n u_m^{-1}
 t_m^{-1}$ be the subgroup element which runs along the other side of
 the strip.  So $w,u$ end in the same plane and both cross the same
 last strip. Now let $w=w'r_nw_n, u=u't_mu_m$. Then $w'y=u'$ are
 stable letter reduced and $w'y$ has $k-1$ stable letters so by
 hypothesis $u',w'y$ and $w'$ cross the same succession of strips and
 planes to get to the same last strip, so we are done.
\end{proof}
\noindent
It follows easily that all stable letter reduced paths to a point have
 parallel stable letter structure, but the result above is stronger.

\begin{cor}
All geodesics for a group element $g\in G$ cross the same succession
of strips and planes from $1$ to $g$.
\end{cor}

Thus to keep track of all potentially shorter paths when reading a
word one can restrict one's attention to a single ``branch'' of the
\cg.

Now lets suppose we are building some (finite) machine that can decide
whether or not a given (stable letter reduced) word is geodesic.  We
start at $1$ in the \cg\ and read until we hit the first stable
letter.  This corresponds to crossing the first strip in the \cg.
\begin{defn}[Witness]
We call another word $u$ a {\em witness} for $w$ if at some point
$\overline{u}= \overline{w(t)}$ and $|u| <|w(t)|$.
\end{defn}
 We think of $u$ as witnessing the fact that $w$ is not a geodesic.
By the corollary all potential witnesses for $w$ must cross this first
strip also.

\begin{figure}[ht!]
  \begin{center}
      \includegraphics[width=7cm]{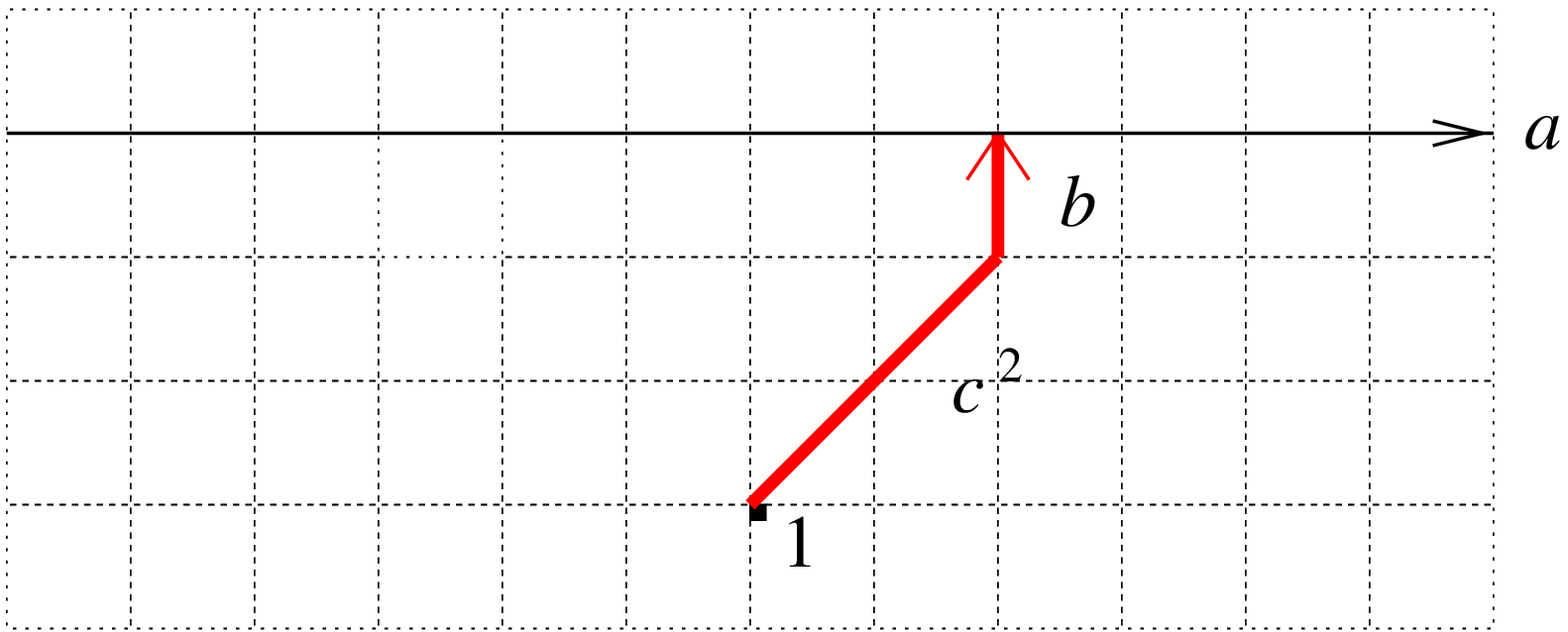}
  \end{center}
  \caption{$w=c^2bs^{-1}w'$ in \bbx}
  \label{fig:2_11extra}
\end{figure}

For example, consider the group \bbx, and the word $w=c^2bs^{-1}w'$
shown in Figure \ref{fig:2_11extra}.
We need to keep track of all the geodesics from $1$ to any ``crossing
 point'' along this strip, since one of these could potentially be a
 witness.  Now it doesn't matter which path one chooses to reach each
 crossing point; we only need to know the length of a geodesic from
 $1$ to this point. In fact, we only need to know the relative
 difference between the length of of geodesic from $1$ to this point,
 and the length of $w$ to its crossing point.

\begin{figure}[ht!]
  \begin{center}
      \includegraphics[width=7cm]{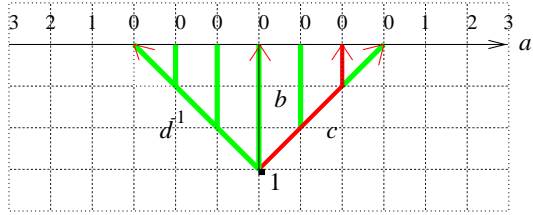}
  \end{center}
  \caption{(All) geodesics to the first strip in \bbx}
  \label{fig:2_11}
\end{figure}

So to keep track of all potential witnesses, we simply record the
relative distances from each crossing point along the strip back to
$1$. For example, in Figure \ref{fig:2_11} we have written these
numbers at each crossing point on the strip.

We are almost ready to make our main definition, but before we do,
recall that we are thinking about building a (finite) machine that
could store this information and so keep track of witnesses as it
reads a word $w$.  So instead of recording the actual distances as we
have written them in Figure \ref{fig:2_11} we instead choose an
orientation of the strip, then record the difference between adjacent
numbers along the strip.
\begin{defn}[Sequence]
A {\em sequence}
is a bi-infinite sequence of numbers ($0,\pm 1$) which correspond to
the difference between the relative distances from adjacent crossing
points along a strip back to 1 in the Cayley graph.
\end{defn}
For example, the sequence for the strip in Figure \ref{fig:2_11} is
$$\ldots , -1,-1, 0,0,0,0,0,0,1,1,\ldots .$$
We write this more neatly as 
$(-1)(0)^6(1)$, where the terminals are always understood to be infinite.

\section{Initial sequences and patterns for \bbx\ and \wx}
\begin{figure}[ht!]
  \begin{center}
  \subfigure{\includegraphics[width=3cm]{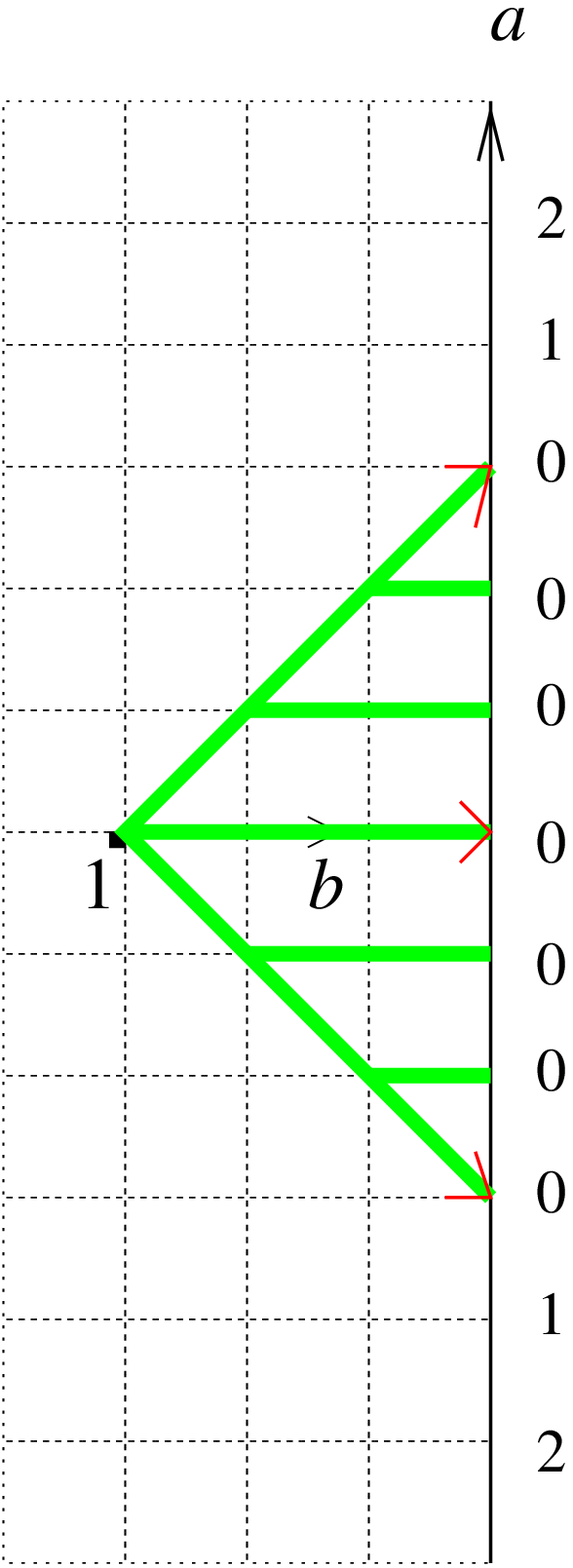}}
\qquad \qquad
 \subfigure{\includegraphics[width=5cm]{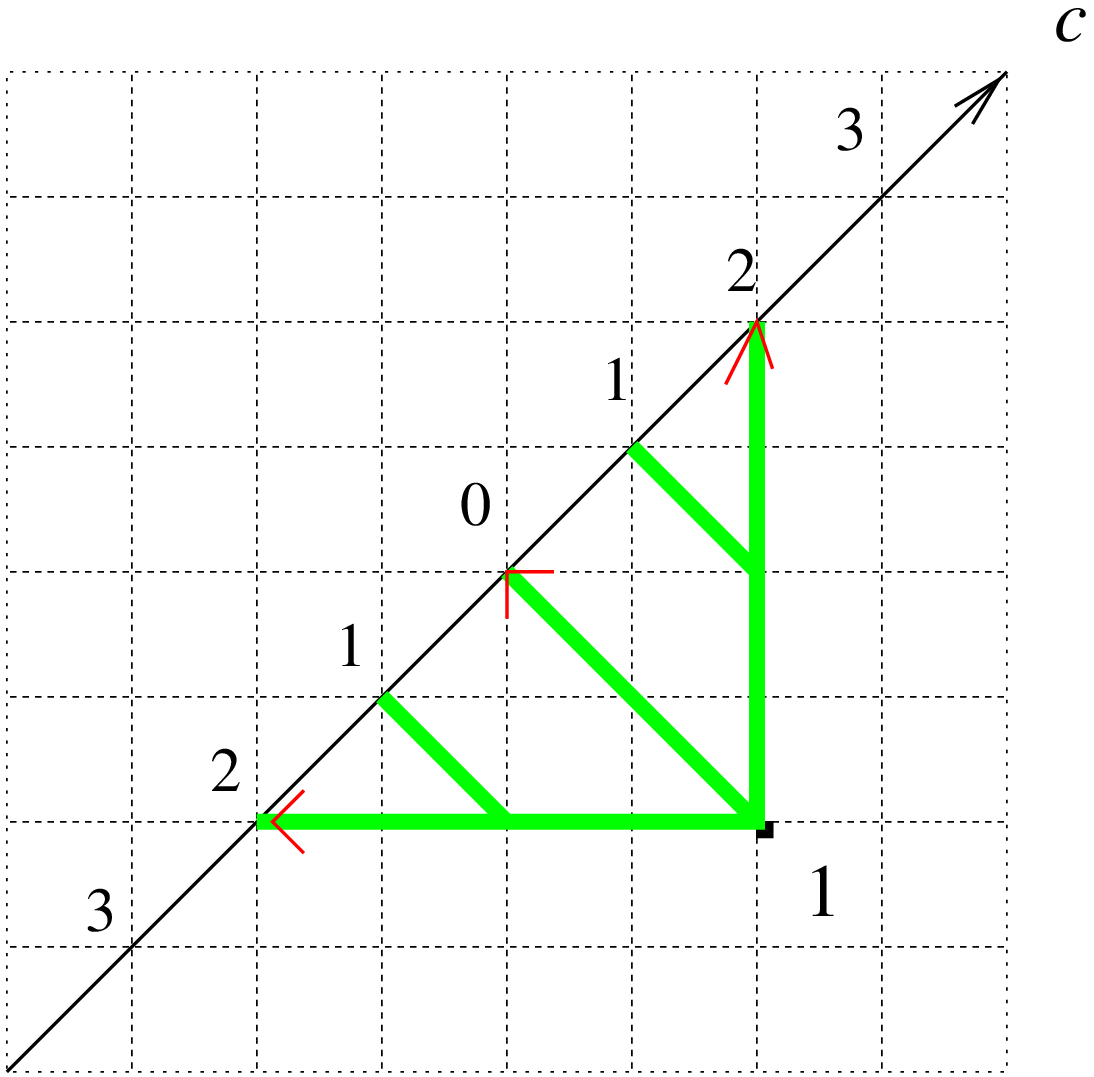}}
\qquad \qquad
 \subfigure{\includegraphics[width=5cm]{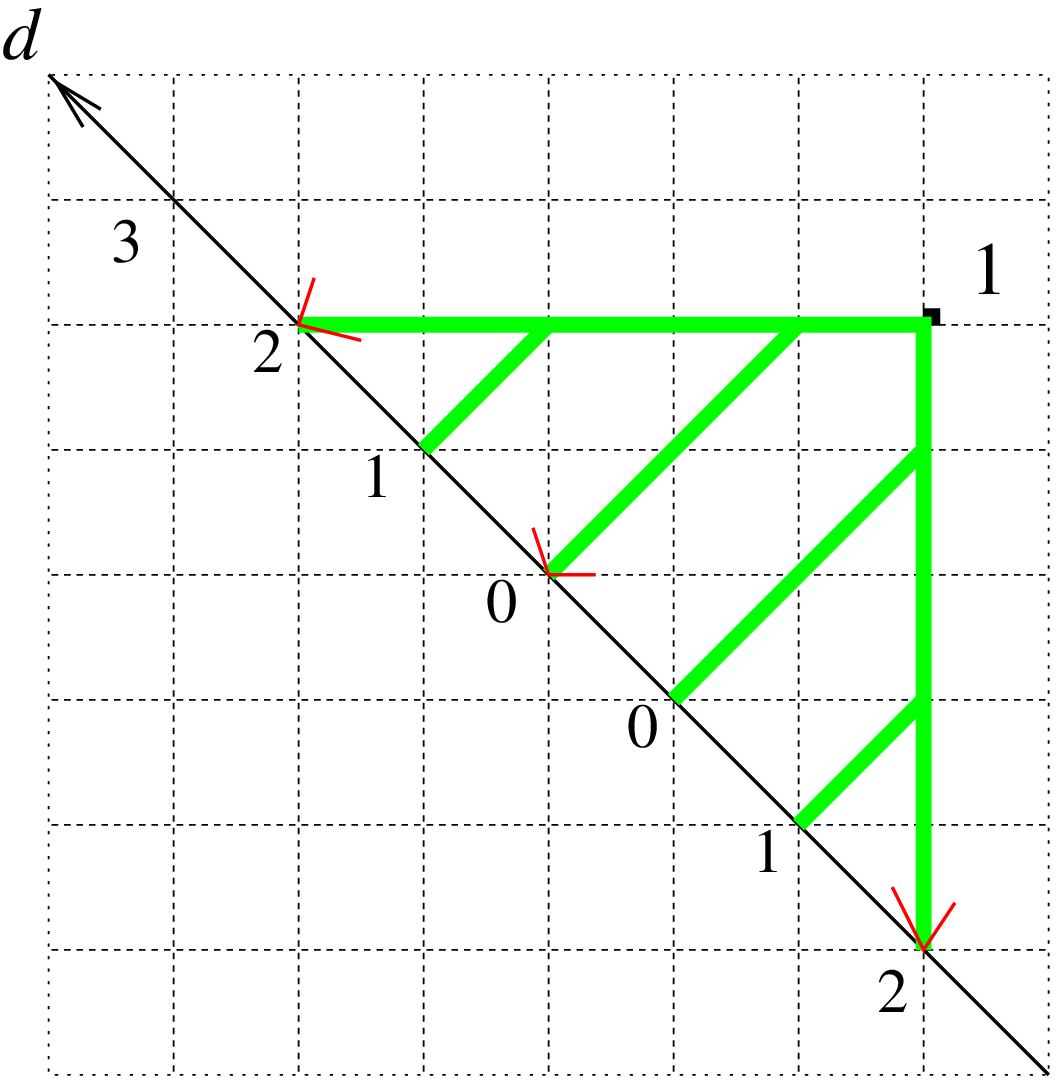}}
  \end{center}
  \caption{Initial sequences for \bbx}
  \label{fig:2_13}
\end{figure}

The diversity of sequences that can arise in a generic multiple \hnn\
 depends on the presentation.  For the group \bbx\ the first strip
 could be glued to either the subgroup $\langle a \rangle, \langle c
 \rangle$ or $\langle d \rangle$.  Each case is shown in Figure
 \ref{fig:2_13}
and we see that all sequences of the form $(-1)(0)^k(1)$ are possible
for any $k=0,1,2i $ ; $ i \in \N$.
 We say that these sequences all
have (or ``belong to'') the same ``pattern'', and we will define the
idea of a pattern presently.  Remember that to build a (finite)
machine which keeps record of certain (bi-infinite) sequences as it
reads a word would need to be able to record arbitrarily large values
of $k$. So instead of recording sequences we could try to record the
patterns.

\begin{defn}[Pattern]
A {\em pattern} 
is an expression of the form

\noindent $(-1)(p_1)(p_2)\ldots (p_k)(1)$ where the $p_i$ are finite
words in $\{0,\pm 1\}$, the $(p_i)$ terms represent any finite number
of repeats of $p_i$, and the terminals are infinite.  A pattern {\em
occurs} in a \cg\ if there is some sequence realizing this pattern
with certain (positive) exponents on each internal term.  The set of
all patterns that can occur for a group presentation will be called
its {\em patterns theory}.
\end{defn}

For \bbx\ each initial sequence has the pattern $(-1)(0)(1)$, as noted
above.

\begin{figure}[ht!]
  \begin{center}
         \subfigure{\includegraphics[height=6.5cm]{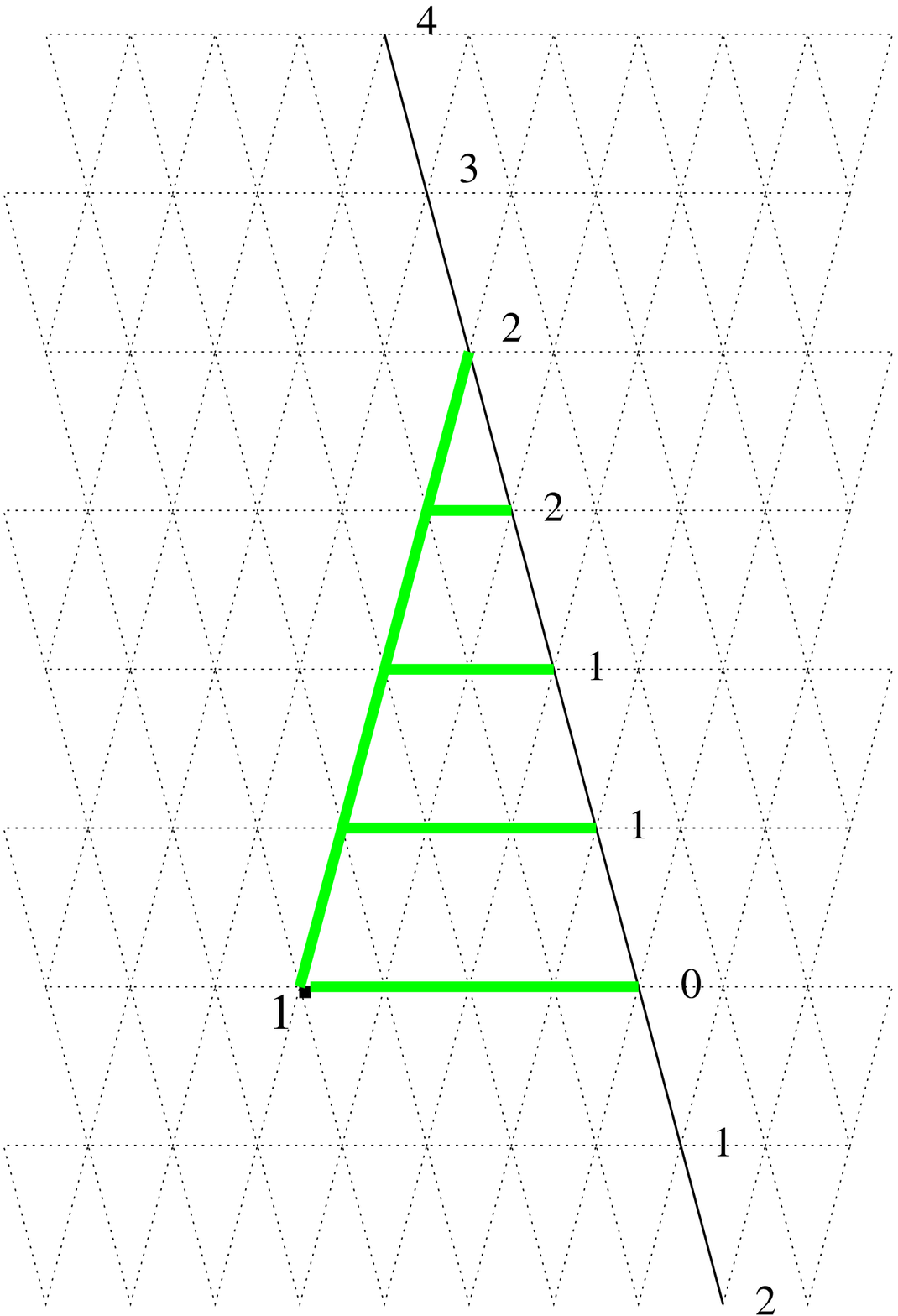}}
         \qquad \qquad
         \subfigure{\includegraphics[width=6.5cm]{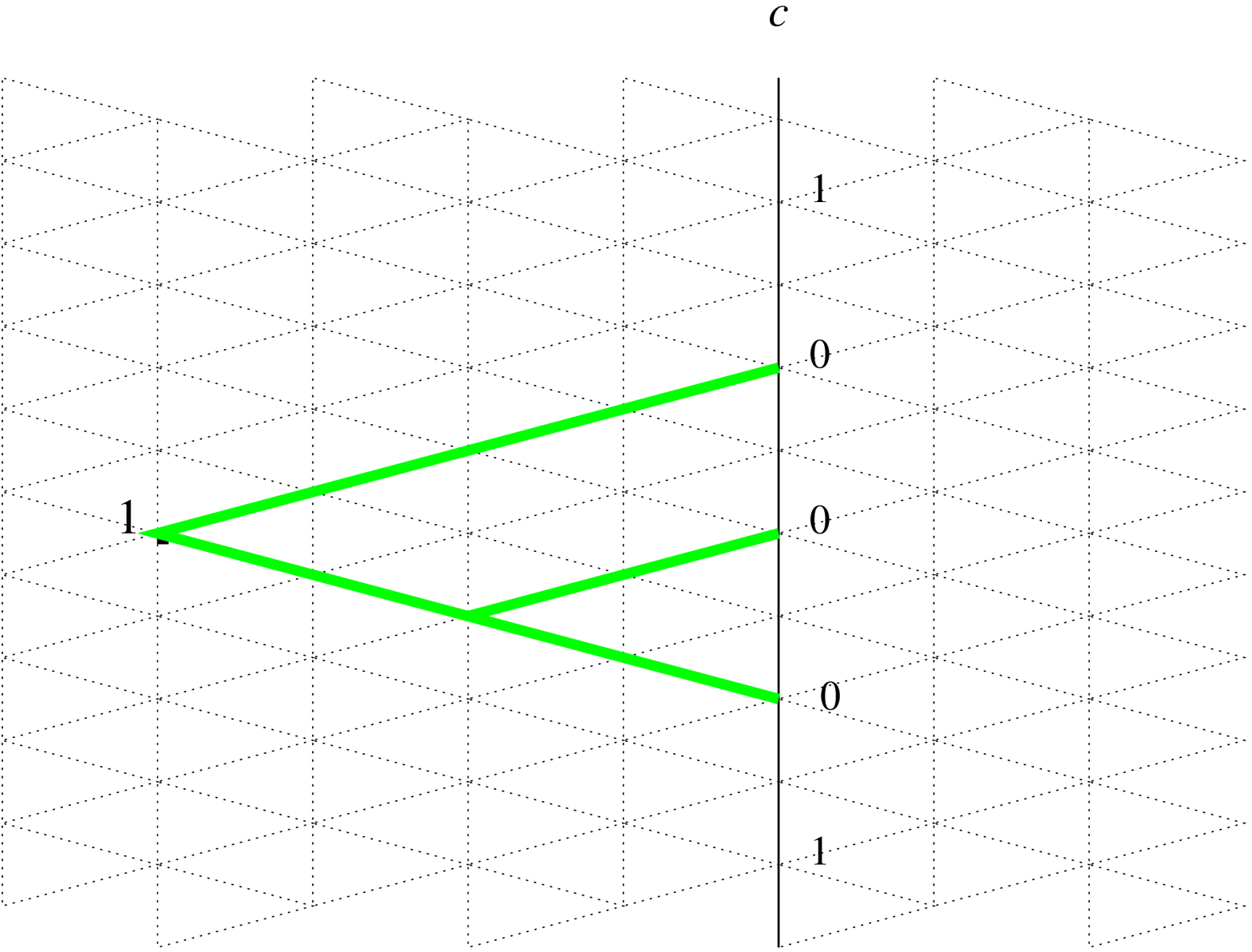}}
         \qquad \qquad
         \subfigure{\includegraphics[width=6.5cm]{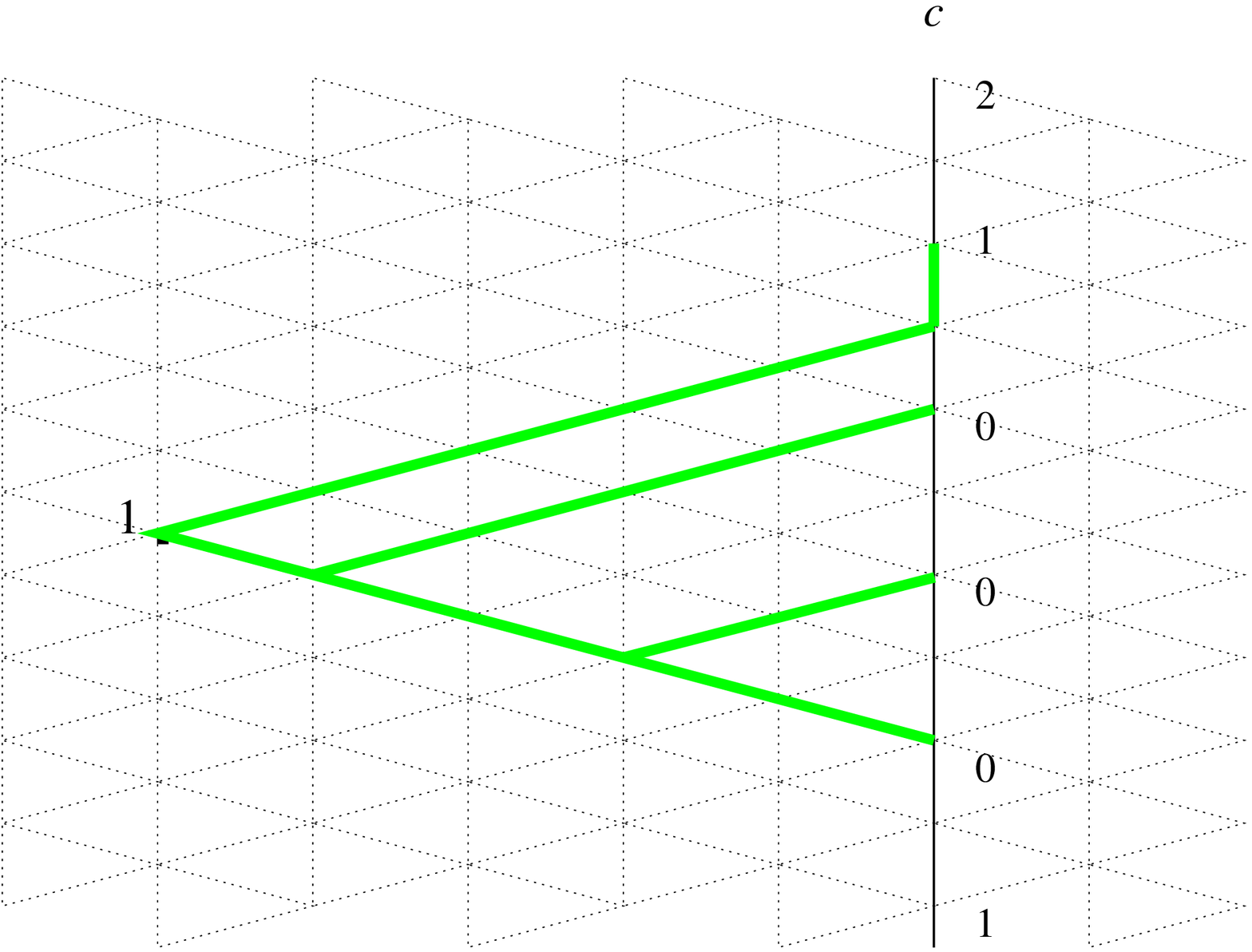}}
  \end{center}
  \caption{Initial sequences for \wx.}
  \label{fig:3_05}
\end{figure}

For the group \wx\ the first strip could be glued to either the
 subgroup $\langle a \rangle, \langle b \rangle$ or $\langle d
 \rangle$.  Each case is shown in Figure \ref{fig:3_05} and we see
 that the possible initial sequences are $(-1)(0)^k(1)$ and
 $(-1)(10)^k(1)$ for any non-negative integer $k$.
\noindent Thus \wx\ has two distinct initial patterns: $(-1)(0)(1)$ and $(-1)(10)(1)$.

So one can see that these two examples have different initial
patterns. Next we will consider what happens on the next strip that a
word crosses, and one might suspect that the new patterns that arise
here will vary in each example.  Therefore to avoid confusion we will
stick with one example and examine the set of all possible sequences
and patterns that can occur.

Also for clarity in the figures we write the actual relative distances
back to $1$, and leave it to the reader to determine the sequence in
each case.

\section{Moves}

We wish to characterize all sequences that occur in the \cg\ of \bbx.
We will prove that only certain patterns occur, that is, all sequences
are of a certain form.

All geodesics start in the plane containing 1.  The strip we use to
exit this plane is called the initial strip.  There are three types of
strips to exit, and all have the same type of pattern, as seen in
Figure \ref{fig:2_13}.

Let an {\em $x$-line} mean a bi-infinite straight path in the \cg\
corresponding to the path $x^i$. The notation $x \ra y$ means that we
enter the plane via a strip glued to an $x$-line and exit the plane by
a strip glued to a $y$-line.

Now there are only four ways to get from one plane to the next:

$\mt 1: a \ra  a, c\ra  c, d\ra  d  $

$\mt 2:  c \ra  a, d\ra  a $

$\mt 3: a \ra  c, a\ra  d $

$\mt 4:  c\ra  d, d\ra  c $

We call these {\em moves} on the patterns, since they take an existing
pattern and make a new one from it.  We wish to write down all
possible patterns that can be generated from the initial pattern by
moves.  Let's start with the initial pattern $(-1)(0)(1)$ and do some
experimentation to see what kind of patterns are possible.

\begin{figure}[ht!]
  \begin{center}
 \subfigure{\includegraphics[width=5cm]{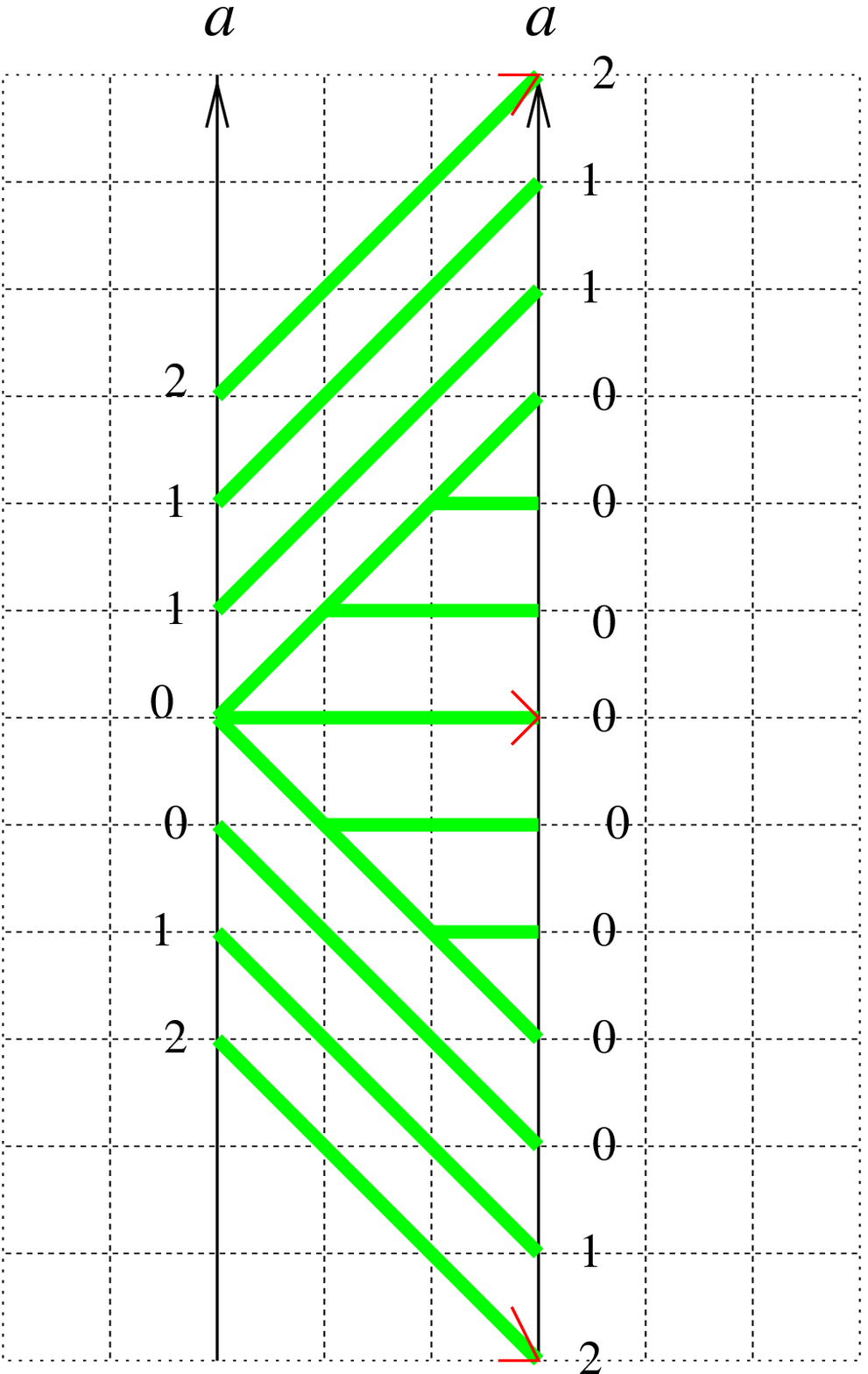}}
\qquad \qquad
 \subfigure{\includegraphics[width=5cm]{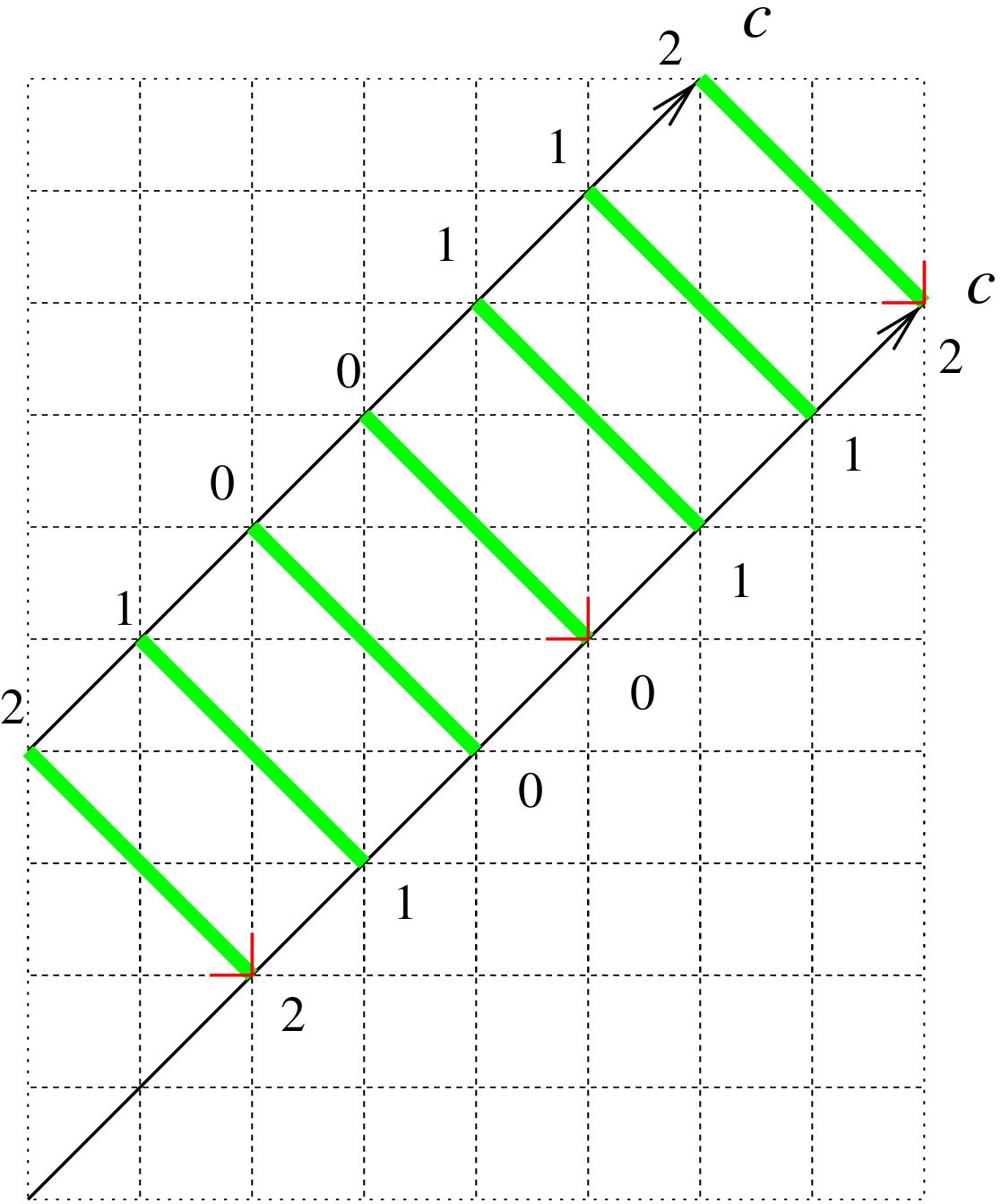}}
\qquad \qquad
 \subfigure{\includegraphics[width=5cm]{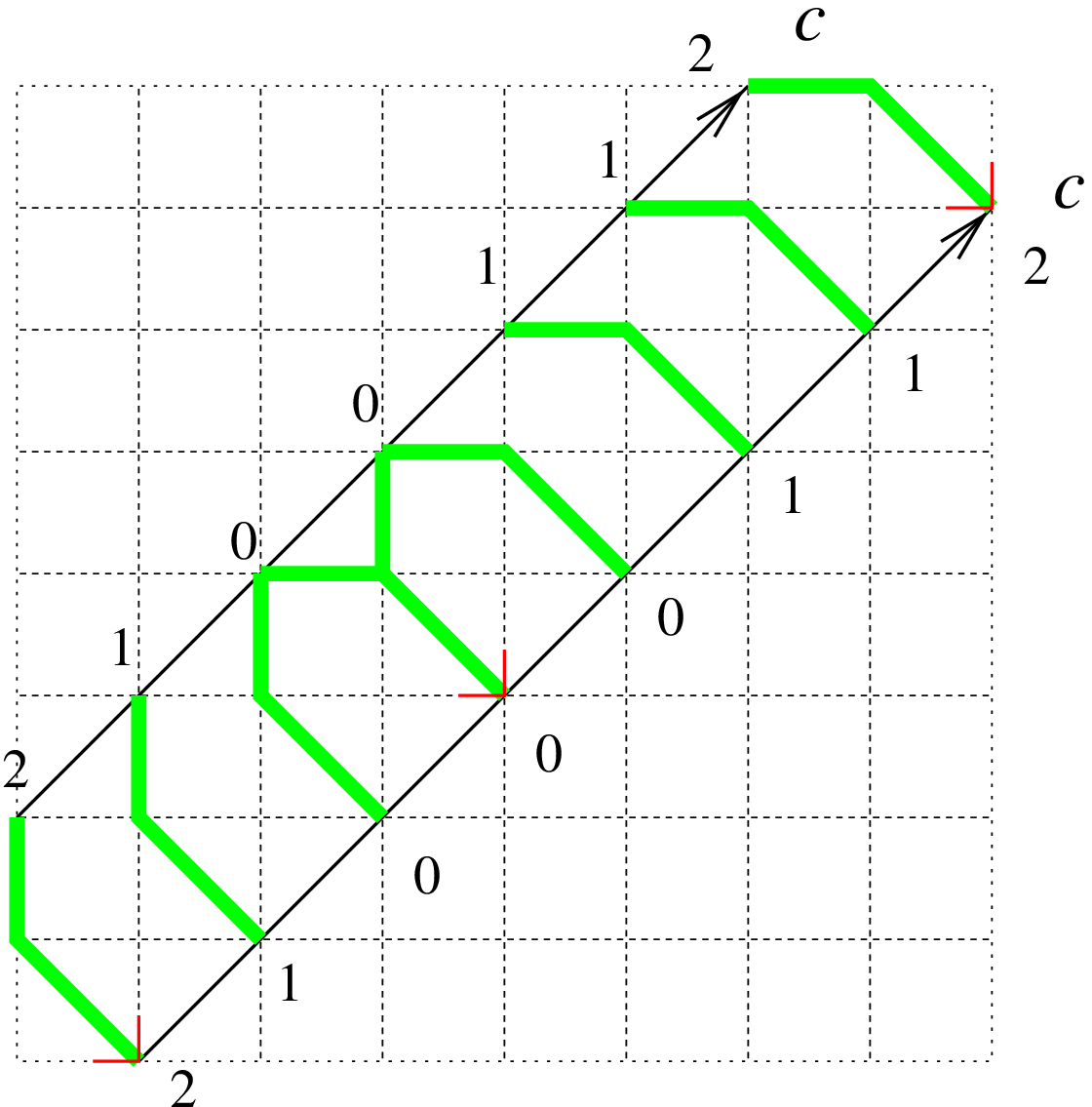}}
  \end{center}
  \caption{``Parallel moves''}
  \label{fig:2_14}
\end{figure}

\noindent
Type $\mt 1$ moves can be referred to as {\em parallel moves}
since they occur when two concurrent strips are parallel on a plane. These
are shown in Figure \ref{fig:2_14}.
In each case we have just ``expanded'' the $(0)$ sub-word of the
initial pattern.  So in this instance Move $\mt 1$ acts as an identity
move since it doesn't introduce any new structure. This may not be
true when we have a more intricate pattern on the enter strip.

\begin{figure}[ht!]
  \begin{center}
      \includegraphics[width=9cm]{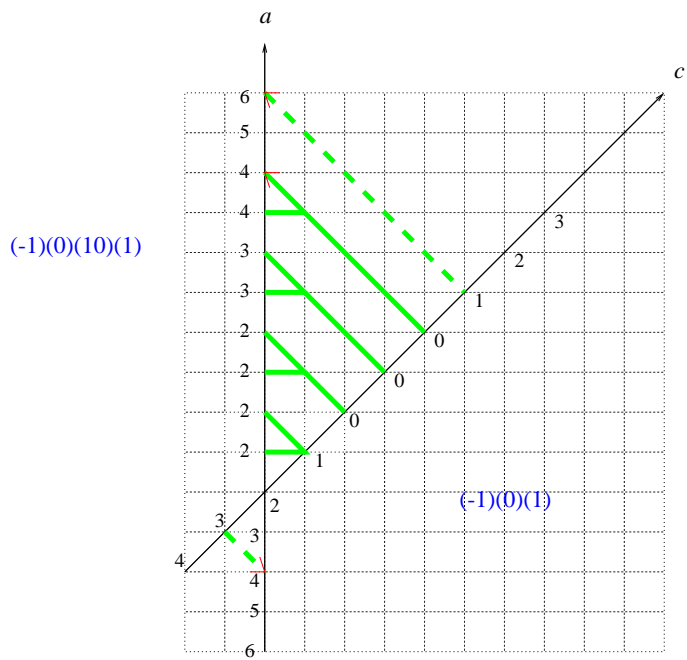}
  \end{center}
  \caption{Move $\mt 2$ gives $(-1)(0)(10)(1)$}
  \label{fig:2_15}
\end{figure}

The type $\mt 2$ move $c/d \ra a$ gives a new pattern as shown in
Figure \ref{fig:2_15}.

\begin{figure}[ht!]
  \begin{center}
      \includegraphics[width=9cm]{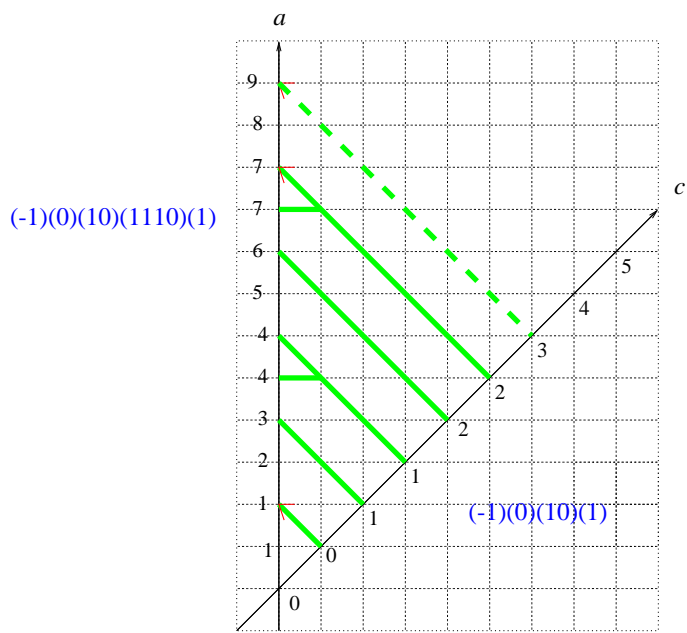}
  \end{center}
  \caption{Another move $\mt 2$ gives $(-1)(0)(10)(1110)(1)$}
  \label{fig:2_16}
\end{figure}

\begin{figure}[ht!]
  \begin{center}
      \includegraphics[width=10cm]{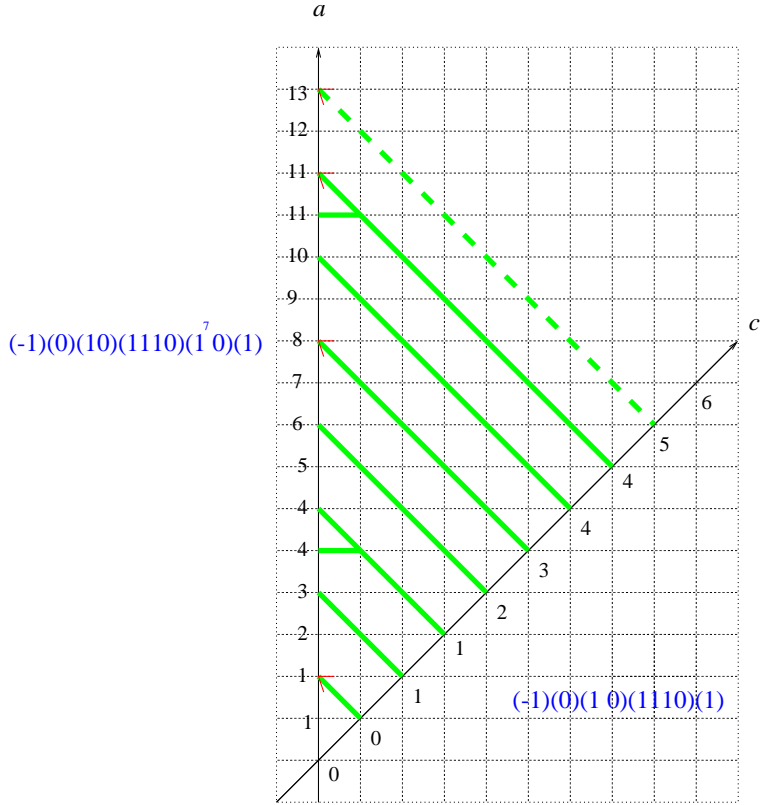}
  \end{center}
  \caption{Another move $\mt 2$ gives $(-1)(0)(10)(1110)(1^7 0)(1)$}
  \label{fig:2_17}
\end{figure}

\noindent
Applying Move $\mt 1$ again just expands the $(0)$ so preserves it.
Applying Move $\mt 2$ again to this gives Figure \ref{fig:2_16}.
Applying Move $\mt 2$ again gives Figure \ref{fig:2_17}.  Note that
the power of seven in the pattern $(1^7 0)$ is not arbitrary in this
case.

Each iteration gives a new term between the two most recently
introduced terms, that is, some kind of ``rewrite'' of the most
recently introduced term.  The reader is encouraged to keep going at
this point, and see what weird and wonderful patterns can be
generated.  What happens when you apply $a \ra c/d,c/d \ra d/c$ and
the parallel moves?  How many jumps of 1 before a 0 are possible?

\begin{figure}[ht!]
  \begin{center}
      \includegraphics[width=7cm]{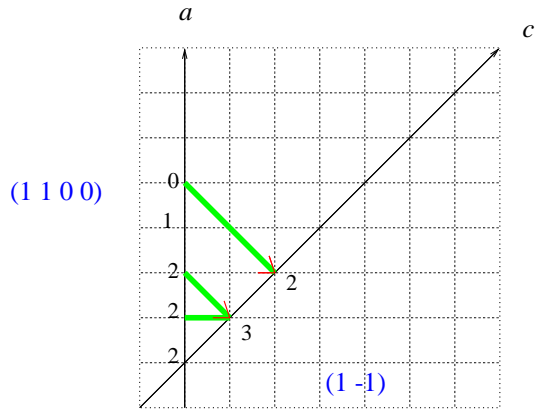}
  \end{center}
  \caption{A potentially ``bad'' pattern}
  \label{fig:2_18}
\end{figure}

We will now introduce a further abstraction of sequences and patterns.
The pattern $$(-1)(0)(10)(1110)(1^7 0)(1)$$ above can be said to
belong to the set of patterns of the form $$(-1)(0)(1,0)(1).$$ where
$(1,0)$ means any mixture of 0's and 1's.
This notation contains less information about the pattern's structure,
but will be useful below.  Note that applying a move $\mt 3$ of the
form $a\ra c$ to this abbreviated pattern potentially can give
$(-1)(0)(1,0,-1)(1)$ as shown in Figure \ref{fig:2_18}, but in the old
notation this would not occur.

\section{Moves as rewriting rules.}

We start with a conjecture about what patterns are allowed.  Using
this we can describe the moves more effectively, in terms of rewrite
rules.
\begin{conj}
All patterns for \bbx\ are of the form
$$(-1)(0,-1)(0)(1, 0)(1).$$
\end{conj} 

 Now assume the pattern $P=(-1)w_1(0)w_2(1)$ is of the conjectured
form, so $w_1,w_2$ are words in $\{-1,0\},\{0,1\}$ respectively.
Define $T=(-1)(0)(1)$ to be the {\em trivial pattern}, which is the
pattern on the first strip for any geodesic.

\subsubsection*{Move $\mt 1$:}
Parallel moves.
\begin{figure}[ht!]
  \begin{center}
 \subfigure[$a \ra a$]{\includegraphics[width=5cm]{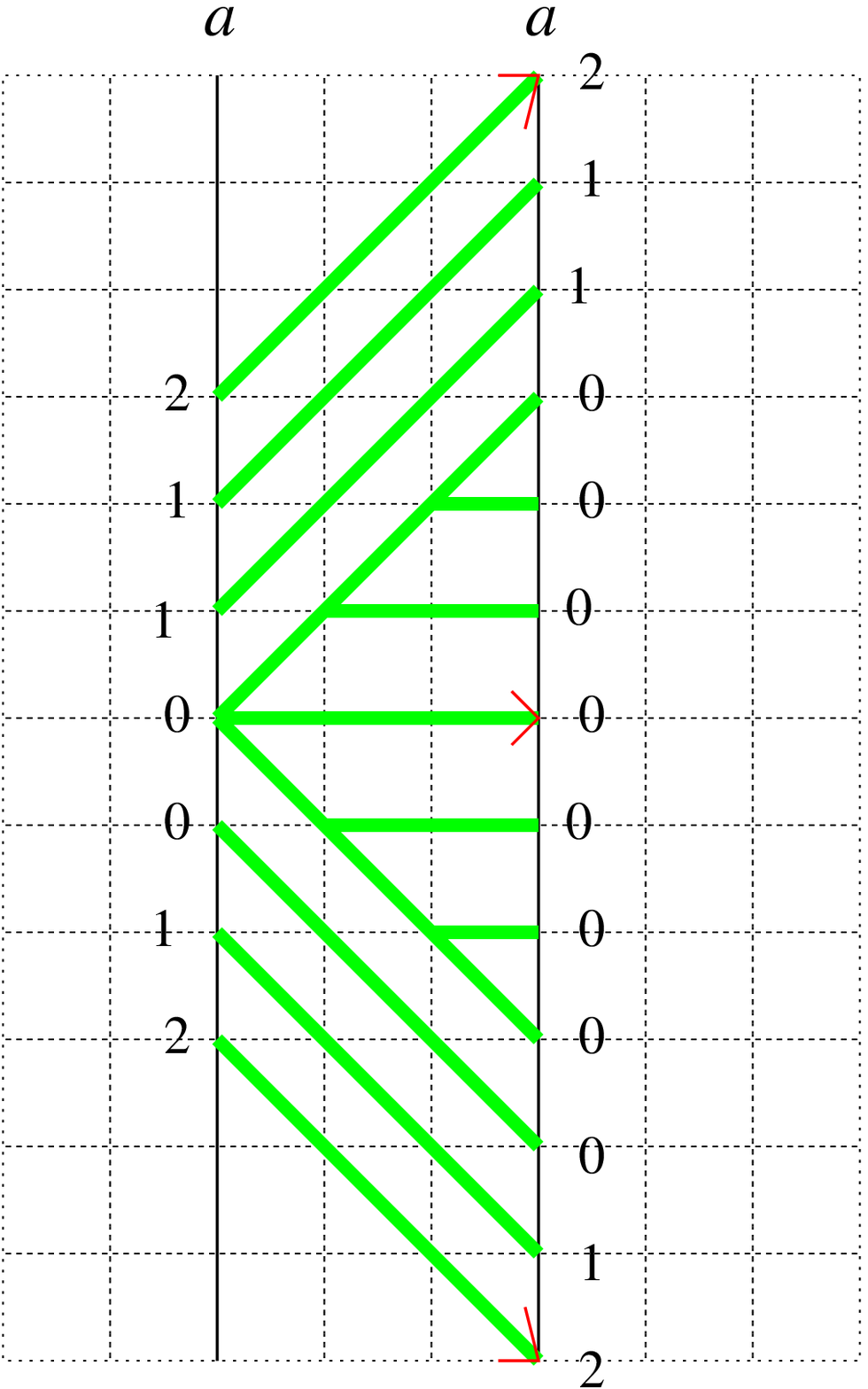}}
\qquad \qquad
  \subfigure[$c / d \ra c / d$]{\includegraphics[width=5cm]{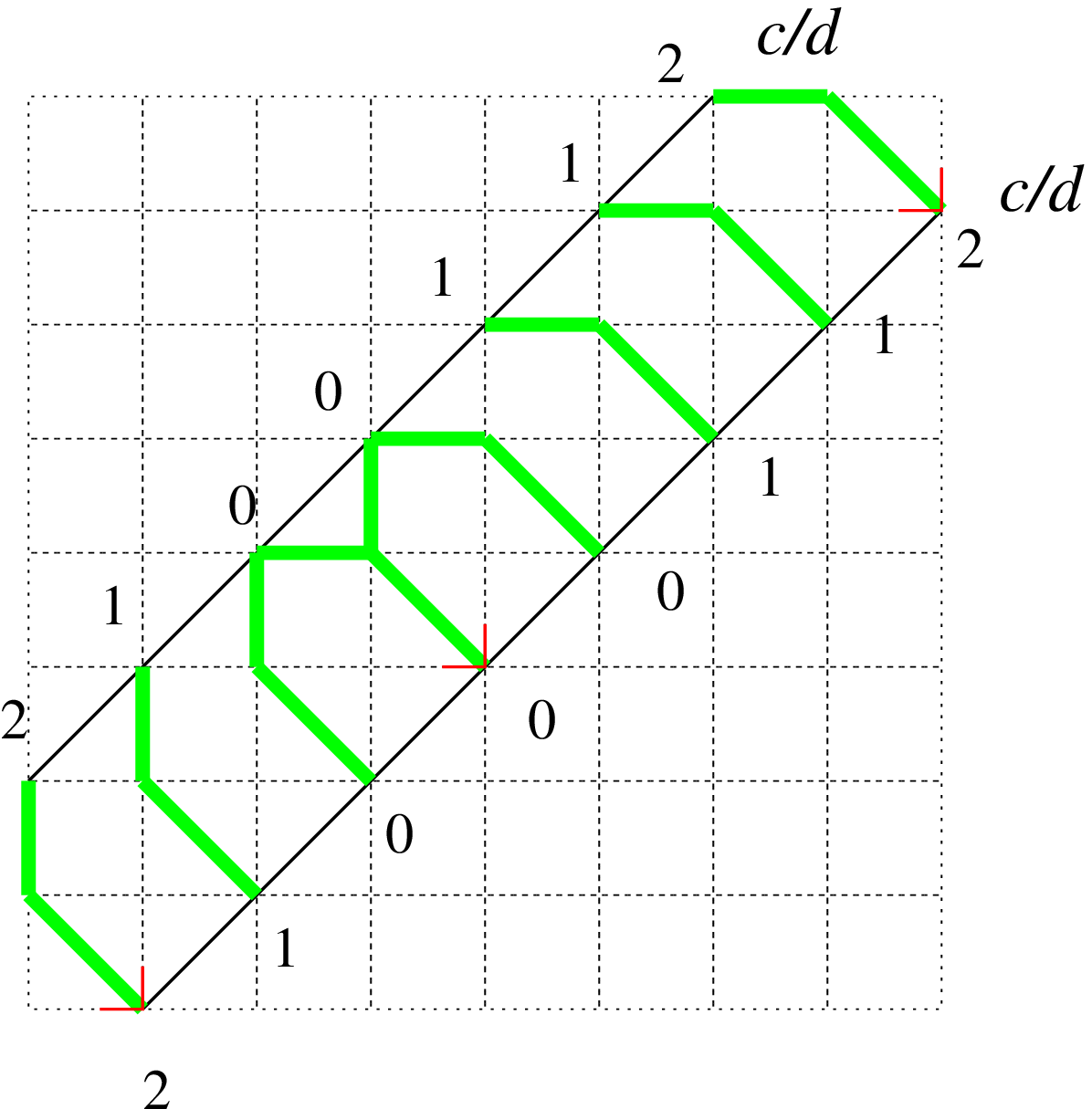}}
  \end{center}
  \caption{Move $\mt 1$}
\label{fig:2_19}
\end{figure}

\noindent
Each of these just increases the number of repeats of $(0)$, so the
pattern is preserved.  Define Move $\mt 1$ to be the {\em identity
move} on the set of patterns, which corresponds to traversing parallel
strips.  So Move $\mt 1$ ``rewrites'' any pattern $P=(-1)w_1(0)w_2(1)$
as itself, since our notation means there can be any finite number of
zeroes in the middle parentheses.

We define a {\em sub-pattern} of $P=(-1)w(1)$ to be a pattern of the
form $(-1)w'(1)$, and a {\em sub-sequence} of $K=(-1)w(1)$ to be a
sequence of the form $(-1)w'(1)$, where $w'$ is a sub-word of $w$.

\subsubsection*{Move $\mt{2}$:} $c / d \ra a$

\noindent
Let $P=(-1)w_1(0)w_2w_3(1)$ be any pattern of the conjectured form.
\begin{figure}[ht!]
  \begin{center} \includegraphics[width=7cm]{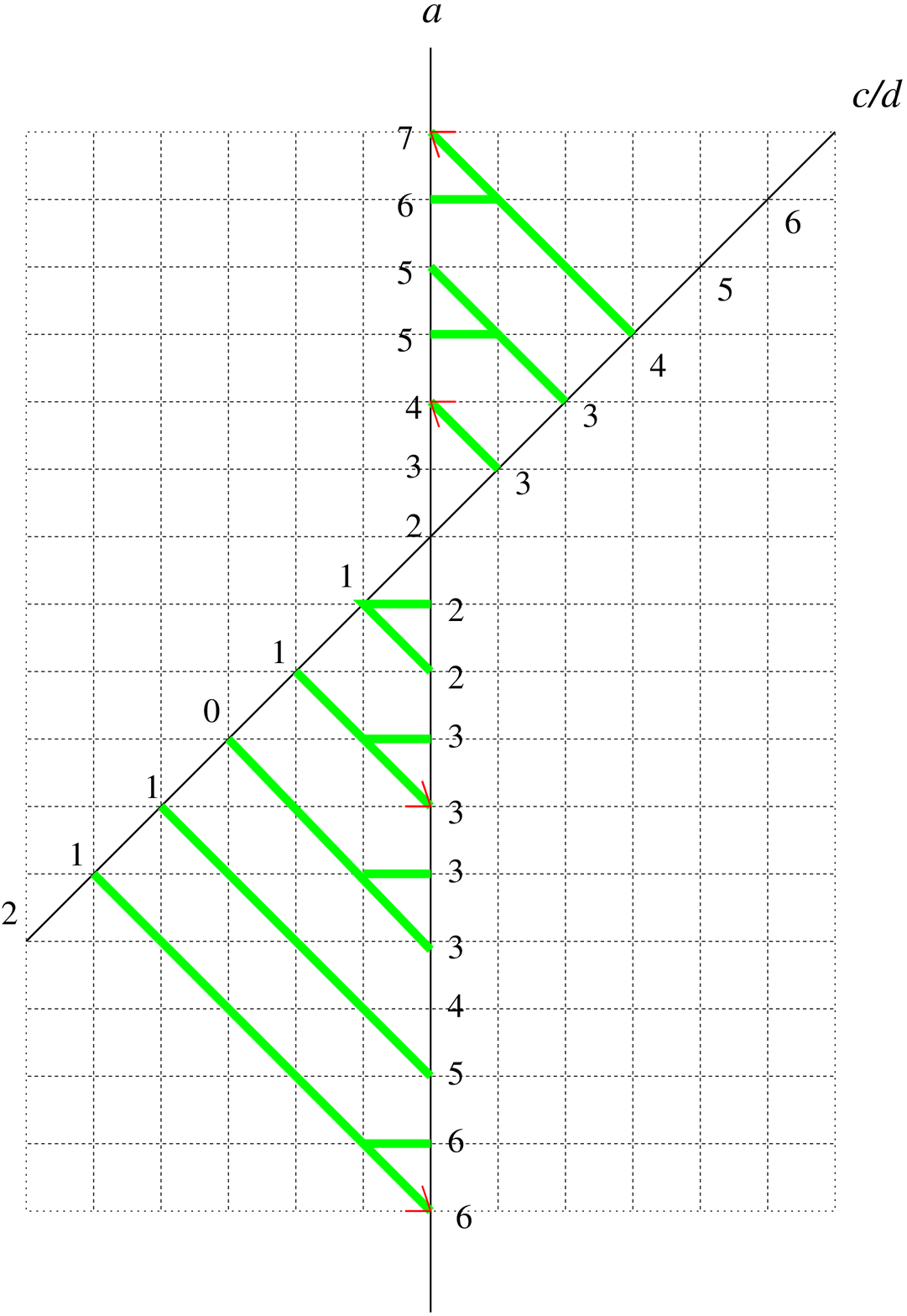}
      \end{center} \caption{Move $\mt{2}$: $c / d \ra a$}
      \label{fig:2_20}
\end{figure}

Rewrite $P$ using the following rules:
\[ w_1(0)w_2 : \left\{ \begin{array}
{r@{\quad \ra \quad}l}
-1 &  -1\; -1 \\ 
0 & 0 \; -1 \\
1 & 0 \; 0 \\
\end{array} \right. \]
\[ w_3 : \left\{ \begin{array}
{r@{\quad \ra \quad}l}
0 & 1 \; 0 \\
1 & 1 \; 1 \\
\end{array} \right. \]
The subdivision of $w_2w_3$ corresponds to the point at which the two
strips intersect on the plane, as in Figure \ref{fig:2_20}.

\subsubsection*{Move $\mt{0}$:}
As defined a pattern has no orientation, so we define a Move $\mt 0$
which rewrites the pattern in reverse.  This move does not correspond
to any strip crossing, it is merely a change of perspective. It is
easily checked that moving $\mt {02}$ covers the remaining ways to go
from $c/d \ra a$.

\subsubsection*{Move $\mt{3}$}: $a \ra c/d$
\begin{figure}[ht!]
  \begin{center}
      \includegraphics[width=7cm]{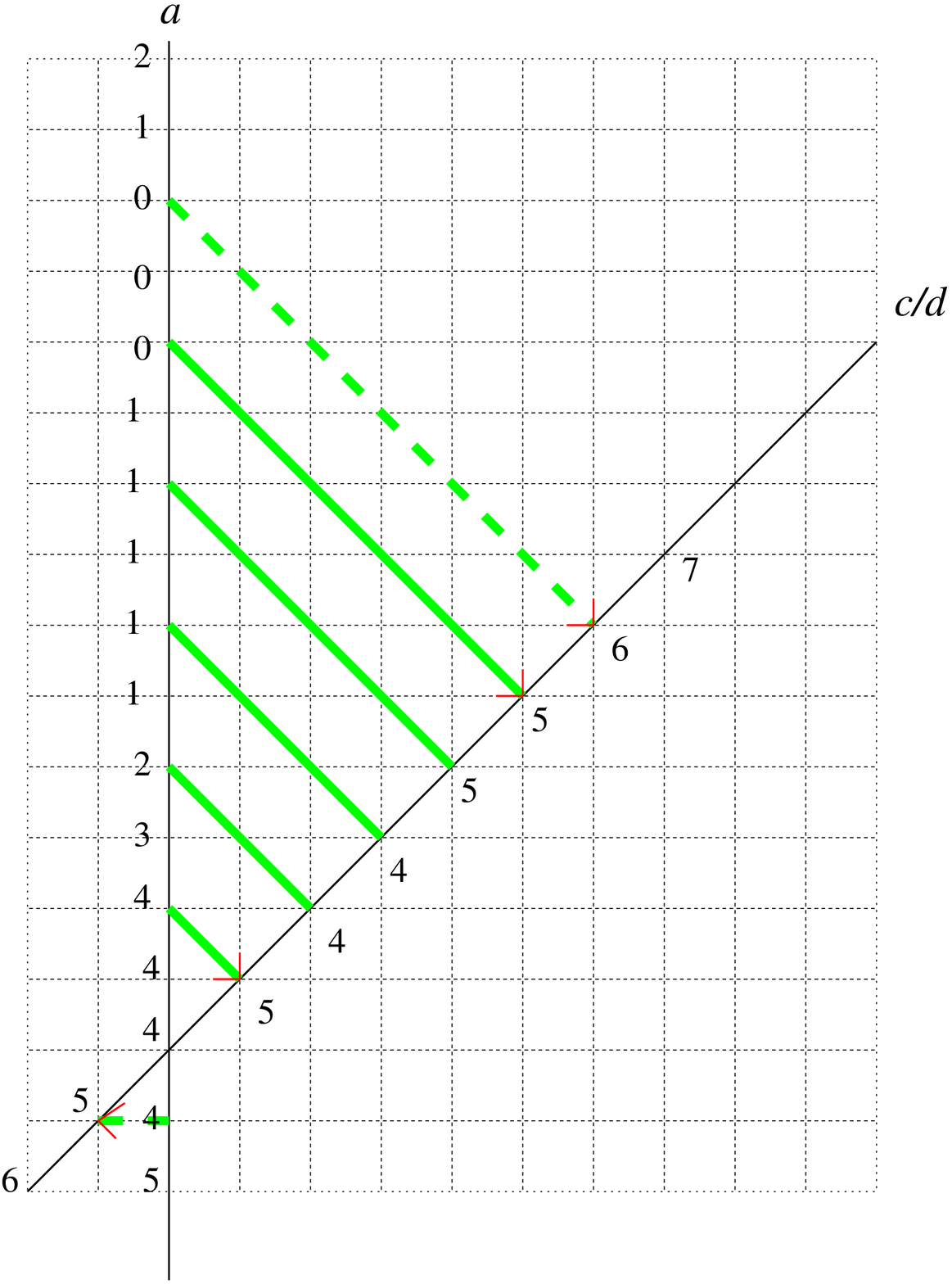}
  \end{center}
  \caption{Move $\mt{3}$: $a \ra c/d$}
  \label{fig:2_21}
\end{figure}

\noindent Take a sub-word $(0)w_1$ of a pattern in conjectured form.
Rewrite:
\[ w_1 : \left\{ \begin{array}
{r@{\quad \ra \quad}l}
0\; 0 & -1 \\
0 \; 1, 1\; 0 & 0\\
1 \; 1 & 1\\
\end{array} \right. \]
The subword $w_1$ in this case corresponds to the part of the pattern
from $(0)$ to the point t which the two strips intersect, as in Figure
\ref{fig:2_21}. Anything ``below'' this does not contribute to the new
pattern.

Notice that this move could potentially give a pattern not in the conjectured form.
For instance if we had a sub-pattern $1100$ which is of conjectured form 
we have seen above this gives $1\;-1$ which is not of conjectured form .

\subsubsection*{Move $\mt{4}$:} $c/d \ra d/c$
\begin{figure}[ht!]
  \begin{center}
      \includegraphics[width=7cm]{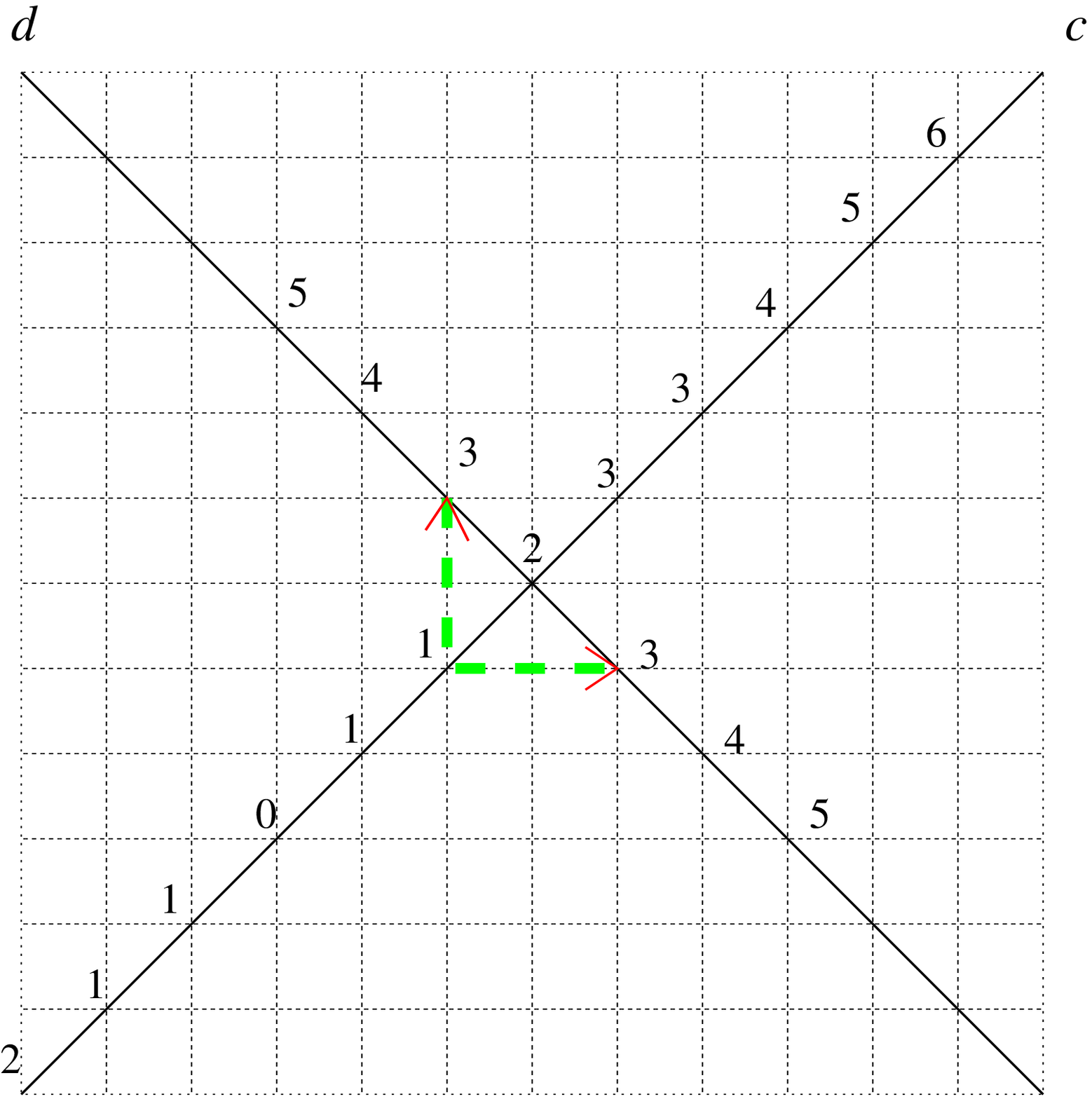}
  \end{center}
  \caption{Move $\mt{4}$: $c/d \ra d/c$}
  \label{fig:2_22}
\end{figure}

\noindent Here we see that regardless of the previous pattern we
always get a pattern $(-1)(1)$ which we can assume is the trivial
pattern.  Therefore only moves $\mt {2,3}$ (and $\mt 0$) can give
non-trivial patterns.

So we now have a way to find every possible pattern for \bbx; start
with the initial (trivial) pattern and perform any number of moves
$\mt {0,2,3}$ in all possible ways, until we obtain a pattern that is
not in conjectured form.  It turns out that all patterns are in
conjectured form, and we will prove this in section \ref{patternsthm}.

\section{Geodesic automatic structures}
Recall that the motivation for this approach was to show that some
geodesic language for \bbx\ can be recognized by a finite machine.

\begin{thm} \label{notreg}
The full language of geodesics for \bbx\ is not regular.
\end{thm}
\begin{proof}
We can find a sequence of the form $(-1)(1^{(2^n-1)}0)(0)$ for
 arbitrarily large values of $n$, by performing some number of move
 $\mt 2$ on the trivial sequence. We show the first three iterations
 of this in Figure \ref{fig:extra1}.  The three segments shown here
 are concatenated together in the \cg.
\begin{figure}[ht!]
  \begin{center}
    \subfigure{\includegraphics[width=8cm]{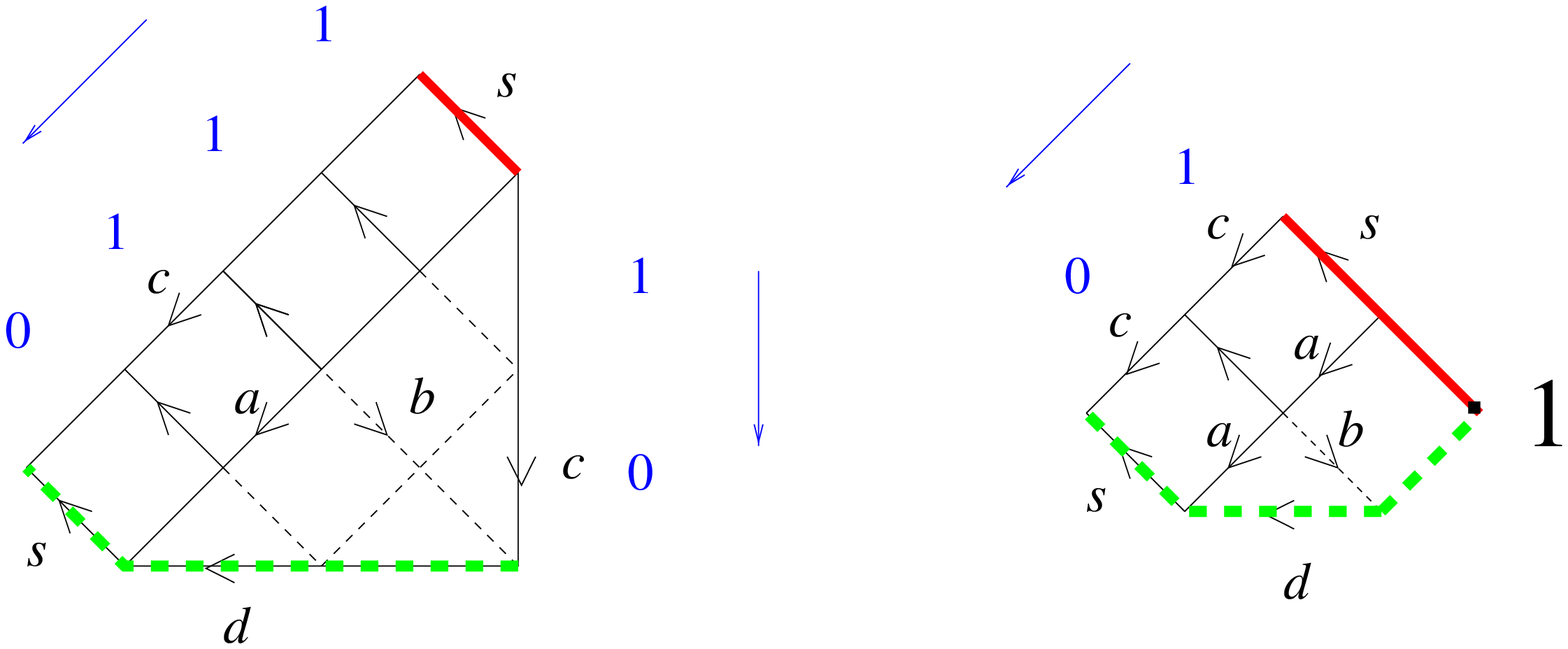}}
\qquad \qquad
    \subfigure{\includegraphics[width=6cm]{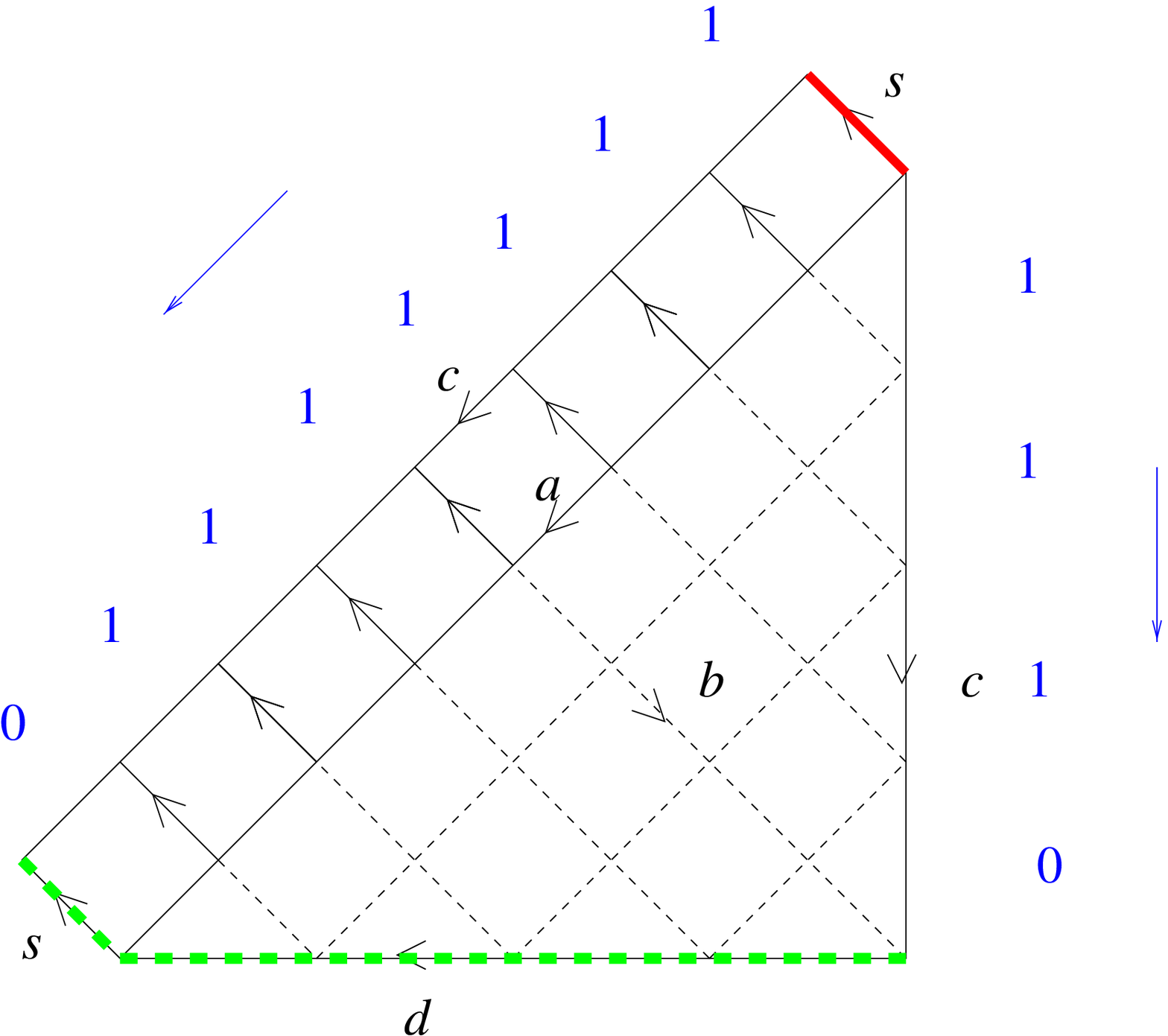}} 
     \end{center}
  \caption{Finding a sequence  $1^k0$}
  \label{fig:extra1}
\end{figure}

This means there is a word $g=b^{-1}s^n$ which is shown in bold in
 Figure \ref{fig:extra1}.  Now $gc^{(2^n-1)}$ is geodesic by the
 labeling of the sequence, and $gc^{2^n}$ is not.  Then by the Pumping
 Lemma the language containing this (geodesic) word cannot be regular.
\end{proof}

\begin{thm}
For any $k>0$ there exist $g,g'\in G_{1,1}$ such that $w\in X^*$ is
the unique geodesic for $g$ and $w'\in X^*$ is the unique geodesic for
$g'$, $d(g,g')=1$ and $w,w'$ do not $k$-fellow travel.
\end{thm}

\begin{cor} \label{nogeodauto}
There is no geodesic automatic structure for \bbx.
\end{cor}
\begin{proof}

We continue with moves on the sequence used in the preceding proof.
Intuitively we can now ``undo'' the $1^k0$ sequence by Move $\mt 3$'s
to get two geodesics which are the unique geodesics from 1 to their
endpoints, ending an edge apart that don't fellow travel.

The previous argument gives us two paths to a strip with sequence
$(1^{(2^n-1)} 0)$.  We wish to perform a move $\mt 3: c/d \ra a$.
First we do a parallel move to get onto a different branch of the \cg,
as in Figure \ref{fig:extra3}.
\begin{figure}[ht!]
  \begin{center}
      \includegraphics[width=7cm]{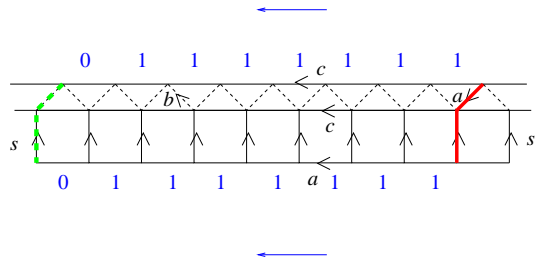}
  \end{center}
  \caption{A parallel move to get on a different branch}
  \label{fig:extra3}
\end{figure}
Then we perform $n$ iterations of move $\mt 3$ as in Figure
\ref{fig:extra4}.  Again these segments are concatenated together
along the appropriate lines in the \cg.
\begin{figure}[ht!]
  \begin{center}
    \subfigure{\includegraphics[width=6cm]{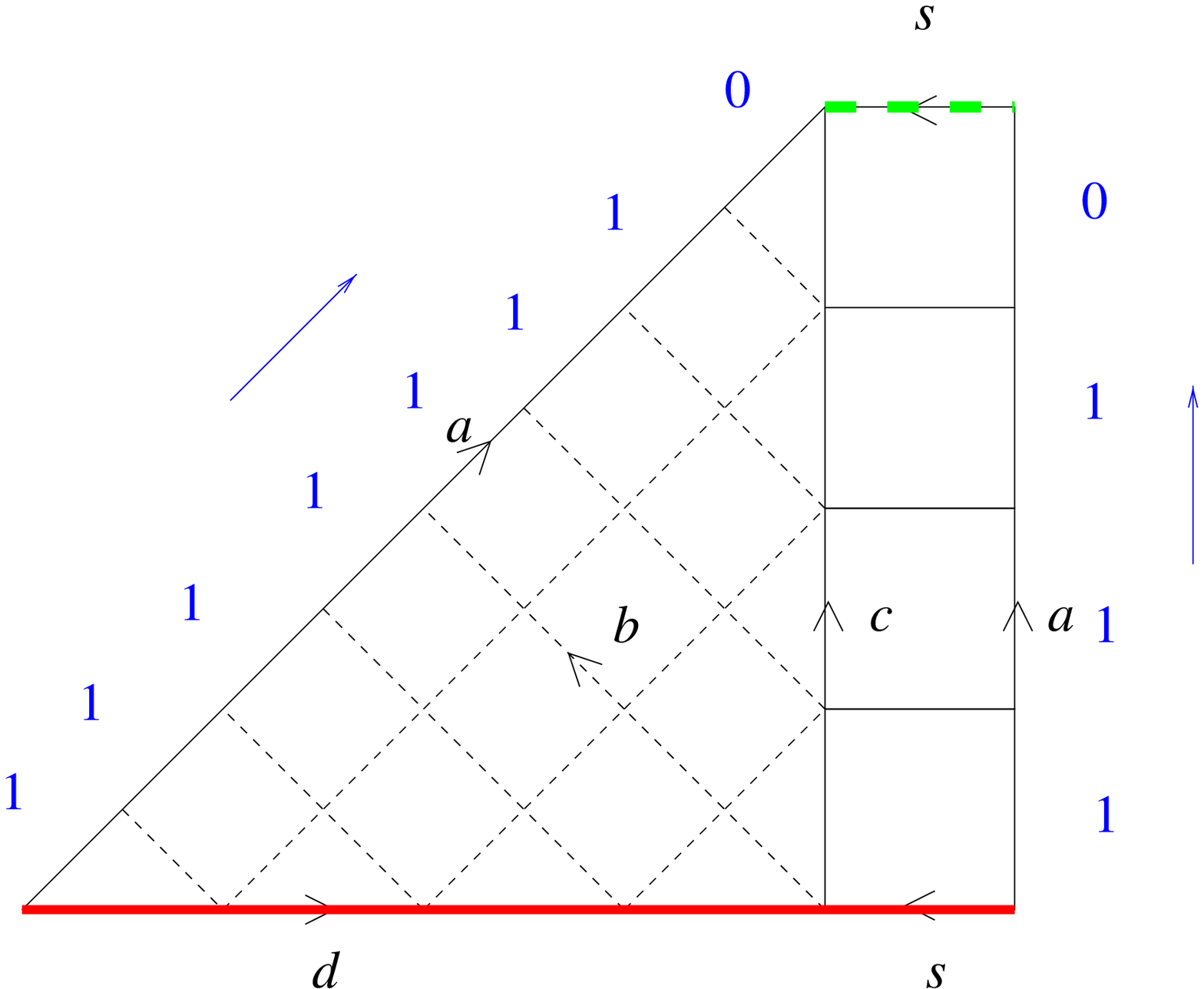}}
 \qquad \qquad
    \subfigure{\includegraphics[width=8cm]{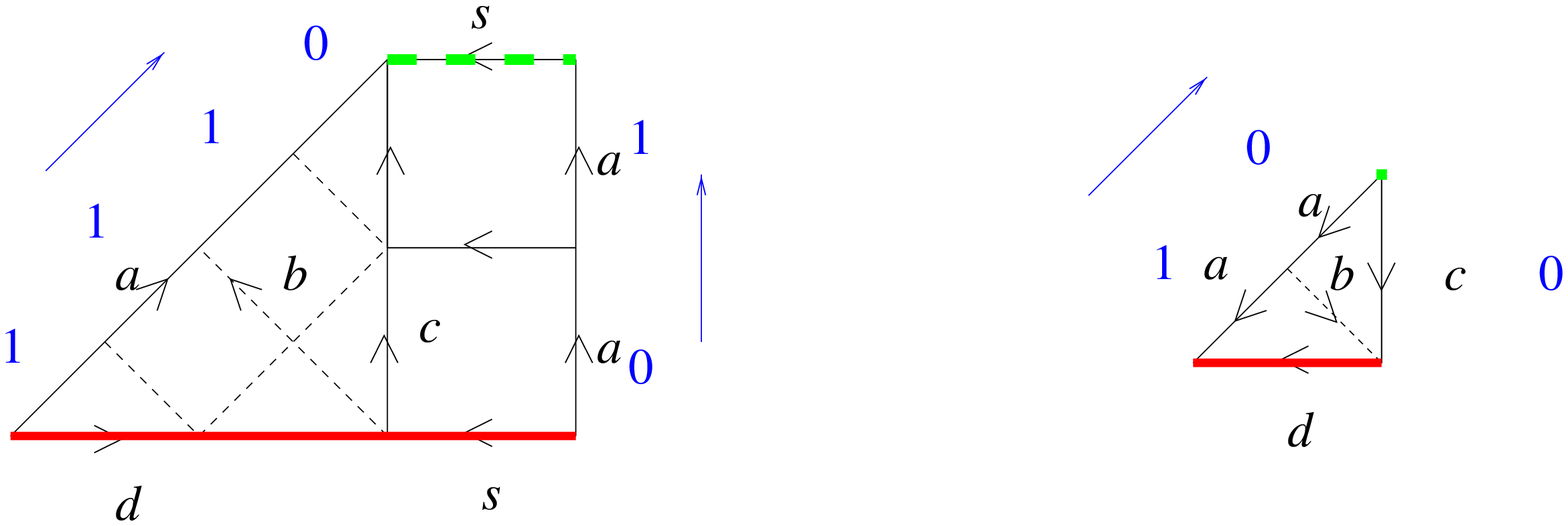}} 
     \end{center}
  \caption{``Undoing'' the pattern}
  \label{fig:extra4}
\end{figure}
Thus we can find two geodesic words
$$w=b^{-1}s^nas^{-1}d^{2^n}s^{-1}d^{2^{(n-1)}}s^{-1}\ldots s^{-1}d^4s^{-1}d^2s^{-1}d$$
$$w'=adsd^2sd^4s \ldots sd^{2^{(n-1)}}sd^{2^n}sas^{-n}, $$
where $w$ is bold and $w'$ is a dashed bold path in the figures.  It
is easily checked, knowing that the sequence on a strip gives all
possible geodesics out to that strip, that $w$ and $w'$ are unique
geodesics to their endpoints which end an edge apart in the \cg, and
fail to $k$-fellow travel for a constant $k$ chosen independently of $
n$.
\end{proof}

\begin{figure}[ht!]
  \begin{center}
      \includegraphics[width=12cm]{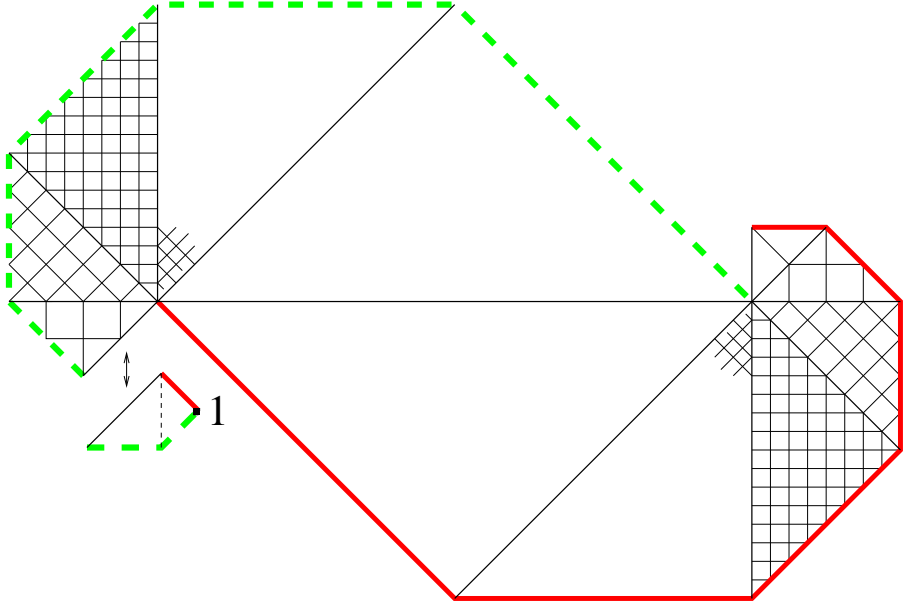}
  \end{center}
  \caption{Unique geodesics that don't fellow travel}
  \label{fig:2_45}
\end{figure}

If we glue the pieces together we get a graphic idea of the geodesic
structure, as in Figure \ref{fig:2_45}.  The parallel move and the
strips have been shrunk to lines for simplicity.

Notice that we have shown \bbx\ has no geodesic automatic language,
 but we have not shown that no geodesic language can be regular, nor
 that some automatic structure containing at least some non-geodesic
 representatives does not exist, so the automaticity of this group
 remains open.

\section{Characterizing patterns} \label{patternsthm}

\begin{thm}[Patterns Theorem]
 \label{patternschar} All patterns for \bbx\ are of the form
$$(-1)(0, -1)(0)(1, 0)(1).$$
\end{thm}

\begin{proof}
For any pattern $P$ we say a succession of moves is {\em efficient} if
the number of moves to get from $T=(-1)(0)(1)$ to a sequence having
pattern $P$ is minimal among all successions of moves that give a
sequence in $P$.  For example, $\mt 0 ^2$ is inefficient.  We will
proceed by induction on the length of efficient successions of moves.

To prove the base step, we know any sequence in $P$ must start at the
trivial pattern.  Applying moves $\mt {0,1,2,3,4}$ to $T$ only $T\mt
2$ is non-trivial, so there is only one pattern $T\mt 2$ with an
efficient succession of one move.

Now suppose $P$ has an efficient succession of $k$ moves.  Suppose $P$
is not of the conjectured form. Then $P$ is non-trivial so an
efficient succession of moves for $P$ starts with a $\mt 2$.  We know
that the only move that can give a bad pattern is move $\mt 3$, thus
the succession ends with a $\mt 3$ or there is a succession with less
moves that gives a pattern not in the conjectured form.  Go to the
last $\mt 2$ in the succession of moves for $P$.  Let $K$ be the
sequence before the last move $\mt 2$.  So the succession of moves
must be $K\mt 2 \ldots \mt 3$.  Since it is efficient, we never get
more than one $\mt 0$ concurrently, so the next moves after $\mt 2$
are either $\mt 3$ or $\mt {03}$.  Let $ K=(-1)w_1(0)w_2w_3(1)$.
$$K\mt 2=(-1)w_1'(0-1)w_2'(0)w_3'(1)$$

\noindent $K\mt {23}:$ Applying move $\mt 3$ only rewrites the $w_3'$.
Let $w_3=x_1\ldots x_n$.  Then $w_3'=y_1\ldots y_n$,
\[ y_i = \left\{ \begin{array}
{r@{\quad : \quad}l}
1\; 0 & x_i=0 \\
1\; 1 & x_i=1 \\
\end{array} \right. \]
$=1a_i$,
\[ a_i = \left\{ \begin{array}
{r@{\quad : \quad}l}
0 & x_i=0 \\
1 & x_i=1 \\
\end{array} \right. \]
So after the move $\mt 3$ we get $(-1)w_3''(1)=(-1)a_1a_2\ldots
a_n(1)=(-1)w_3(1)$.  Thus if $K\mt {23} $ is not of conjectured form
then neither was $K$, which is contradicts the inductive assumption.

\noindent
$K\mt  {203}: $ This time the $w_1'(0-1)w_2'$ part is the only part rewritten.
Let $$K=(-1)x_1\ldots x_n (0)x_{n+1}\ldots x_{m}w_3(1)$$
$$K\mt 2=(-1-1)y_1\ldots y_n (0-1)y_{n+1}\ldots y_{m}(0)w_3'(1)$$
\[ y_i = \left\{ \begin{array}
{r@{\quad : \quad}l}
-1\; -1 & x_i=-1 \\
0\; -1 & x_i=0 \\
0\; 0 & x_i=1 \\
\end{array} \right. \]
Apply move $\mt  0$:
$$\ldots(0)p_m\ldots p_{n+1}(10)q_n\ldots q_1$$
\[ p_i = \left\{ \begin{array}
{r@{\quad : \quad}l}
1\; 0 & x_i=0 \\
0\; 0 & x_i=1 \\
\end{array} \right. \]
$=a_i0,$
\[ a_i = \left\{ \begin{array}
{r@{\quad : \quad}l}
1 & x_i=0 \\
0 & x_i=1 \\
\end{array} \right. \]
\[ q_i = \left\{ \begin{array}
{r@{\quad : \quad}l}
1\; 1 & x_i=-1 \\
1\; 0 & x_i=0 \\
\end{array} \right. \]
$=1b_i,$
\[ b_i = \left\{ \begin{array}
{r@{\quad : \quad}l}
1 & x_i=-1 \\
0 & x_i=0 \\
\end{array} \right. \]
Now apply move $\mt 3$: There are two choices for pairing, lets write
it out:
$$\ldots 000a_m0\ldots a_{n+1}010\ldots 101b_n\ldots 1b_1 111\ldots$$
Either choice of pairing gives
$$(-1)z_m\ldots z_{n+1}(0)z_n\ldots z_s(1), s\leq m$$
\[ z_i = \left\{ \begin{array}
{r@{\quad : \quad}l}
-1 & a_i=0, x_i=1 \\
 0 & a_i=1, x_i=0 \\
 0 & b_i=0, x_i=0 \\
 1 & b_i=1, x_i=-1 \\
\end{array} \right. \]
\[= \left\{ \begin{array}
{r@{\quad : \quad}l}
-1 &  x_i=1 \\
 0 &  x_i=0 \\
 1 &  x_i=-1 \\
\end{array} \right. \]
Thus $K\mt {203}$ is a subsequence of $K\mt 0$.  So if $K\mt {203}$ is
not of conjectured form then neither is $K\mt 0$, but $K\mt 0$ has a
shorter efficient succession of moves which contradicts the inductive
assumption.
\end{proof}

Thus all sequences are of a reasonable form.  It would be interesting
to characterize the types of sequences that occur. Would they be
recognized by some computing machine?  That is, to which formal
language class might the set of sequences for \bbx\ belong?

\section{Results for \wx}
Recall that the initial patterns for \wx\ are 
$(-1)(0)(1)$ and $(-1)(10)(1)$.
There are 5 types of moves for this example.

$\mt 0 :$ Reverse orientation of pattern

$\mt 1 : a \ra  a, b\ra  b, d\ra  d  $

$\mt 2 :  d \ra  a, d\ra  b $

$\mt 3 :  a \ra  d, b\ra  d $

$\mt 4 :  a\ra  b, b\ra  a $

Note that move \texttt{4} is now non-trivial, in contrast to \bbx.
Following a similar program for this example we can prove a ``patterns
theorem'', that the full language of geodesics is not regular, and
that no geodesic automatic structure exists for \wx.  Details are left
to the enthusiastic reader, and the same results can be found in
\cite{\Ethesis} for a slightly different (weighted) \gset.

\section{Almost convexity}
In this section we prove an \ac ity result for multiple \hnn s which
have ``well behaved'' patterns, that is, which satisfy a condition
like that in Theorem \ref{patternschar}. The proof parallels that
given by the author in \cite{\Enonhopf} for different hypotheses.

\begin{defn}[Almost convex]
$(G,X)$ is {\em almost convex} if there is a constant $C$ such that
every pair of points lying distance at most 2 apart and within
distance $N$ of the identity in $\Gamma(G,X)$ are connected by a path
of length at most $C$ which lies within distance $N$ of the identity.
\end{defn}

\noindent
See \cite{\Cannon} for properties of \ac\ groups.

\begin{defn}[Falsification by fellow traveler property]
$(G,X)$ enjoys the {\em \fftp} if $\exists k$ such that every
non-geodesic word $w\in X^*$ is $k$-fellow traveled by a shorter word
with the same start and end points.
\end{defn}

\noindent
If $(G,X)$ has the \fftp\ then the language of all geodesics on $X$ is
regular \cite{\NSgeomfinite}. It follows that \bbx\ does not enjoy
this property. If $(G,X)$ has the \fftp\ then it is \ac, with constant
$C=3k$. Both properties depend on the choice of \gset.

\begin{defn}[Well behaved]
\label{wellbehaved}
Let $(G,X)$ be a multiple \hnn\ as in Definition \ref{mult} with a
strip equidistant presentation.  We say $(G,X)$ has {\em well behaved
patterns} if each pattern is of the form $(-1,0)(0)(1,0)$.
\end{defn}


\begin{thm}
\label{acthm}
Let $(G,X)$ be a multiple \hnn\ of $(A,Z)$ as in Definition \ref{mult}
with a strip equidistant presentation, such that associated subgroups
are geodesic, $(A,Z)$ enjoys the \fftp\ and $(G,X)$ has well behaved
patterns.  Then $(G,X)$ is \ac.
\end{thm}

\noindent
{\em Remark:} In \cite{\Enonhopf} the well behaved patterns hypothesis
 is absent, and instead we require that associated subgroups are
 totally geodesic. Wise's example satisfies either set of hypotheses,
 while \bbx\ has well behaved patterns but fails to have totally
 geodesic associated subgroups, as noted in Examples
 \ref{tot1},\ref{tot2} above. Neumann and Shapiro prove that an
 abelian group with any finite \gset\ has the \fftp\
 \cite{\NSgeomfinite}, so it follows that both examples are \ac.

\begin{proof}
Let $S(N)$ denote the metric sphere of radius $N$ and $B(N)$ the
 metric ball of radius $N$ in $\Gamma(G,X)$.  Let $g,g'\in S(N)$ with
 $d(g,g')\leq 2$ realized by a path $\gamma$.  Let $w,u$ be geodesic
 words for $g,g'$ respectively.  Since the presentation is strip
 equidistant, $w,u$ are stable letter reduced.  Let $k$ be the \fftp\
 constant for $(A,Z)$ and assume that $k$ is an even integer greater
 than $\max\{|u_{i_j}|:u_{i_j} \mathrm{is} \; \mathrm{a} \;
 \mathrm{generator} \; \mathrm{of}\; U_i \; \forall i \}$.

\subsection*{Case 1} 
$w, u, \gamma$ have no stable letters. The word $w\gamma$ is not
geodesic in $(A, X)$, so applying the \fftp\ we can find a shorter
word $q$ which ends at $g'$ and $k$-fellow travels it.  If $q$ is not
geodesic then we can find a shorter word $y$ which ends at $g'$ and
$k$-fellow travels $q$.  Moreover $y$ is a geodesic so has length $N$.
The path that retraces $w$ back to $w(N-\frac{k}{2})$, then across to
$q(N-\frac{k}{2})$ by a path of length at most $k$, then across to
$y(N-\frac{k}{2})$ by a path of length at most $k$, then along $y$ to
$g'$, stays within $B(N)$ and has length at most $3k$.

\subsection*{Case 2}
 $\gamma$ involves a stable letter.  Then $\gamma =s,sx, xs$ for $s\in
\{s_i^{\pm 1}\}_{i=1}^n$ and $x$ any generator or inverse of a
generator except $s^{-1}$.  By Britton's Lemma $w\gamma u^{-1}$
contains a pinch so either $w$ or $u$ has an $s^{-1}$.  Without loss
of generality assume $w=w_1s^{-1}w_2$.

\subsubsection*{Case 2a}
If  $\gamma =s, sx$ then we have Figure \ref{figA}(a).
\begin{figure}[ht!]
  \begin{center}
         \subfigure[Case 2a]{
	 \includegraphics[width=5cm]{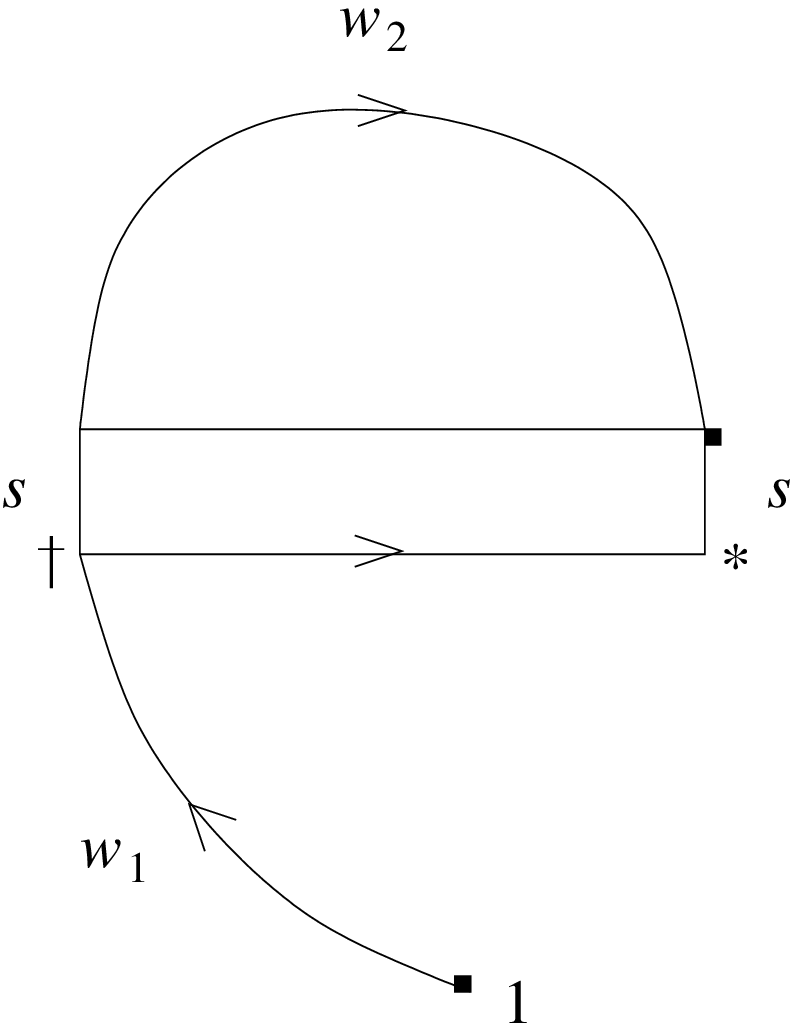}}
	\qquad \qquad
       \subfigure[Case 2b]{
 \includegraphics[width=5cm]{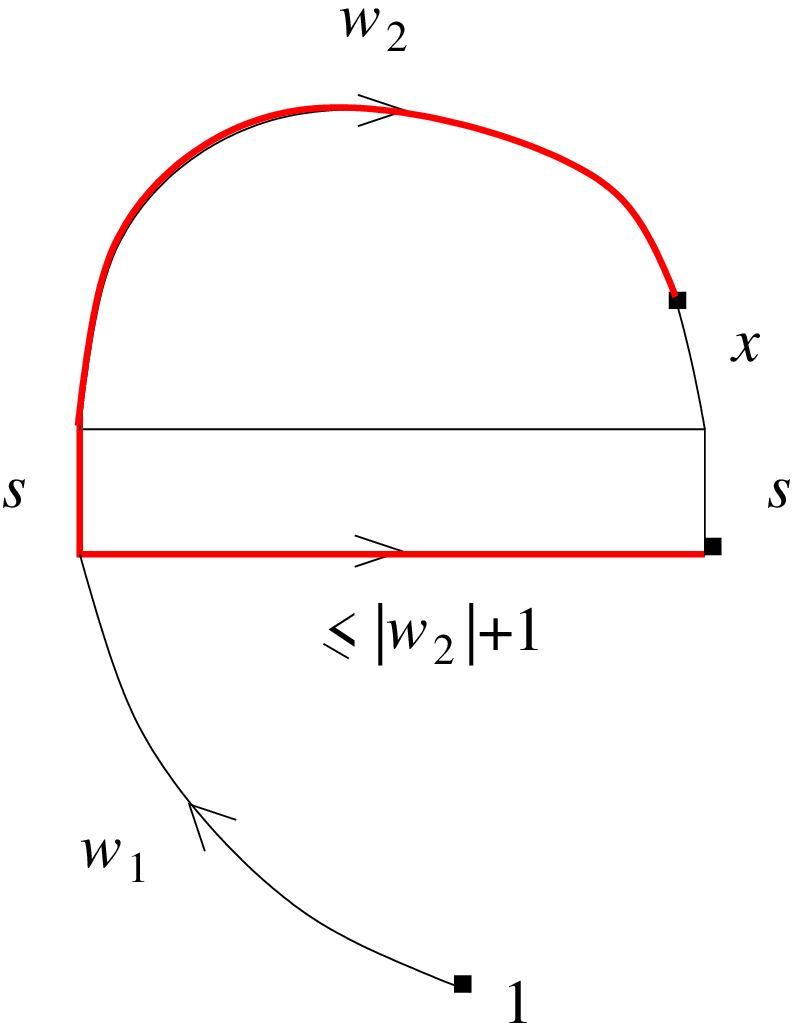}}
	   \end{center}
  \caption{}
 \label{figA} 
\end{figure}
$w_2$ is geodesic and evaluates to an element of an associated
subgroup $U_i$ or $V_i$. 
of $U_i$.  The point $\overline{w_1}=\dag$ lies in $B(N-|w_2|-1)$ and
since associated subgroup words are geodesic, the point $\overline
{ws} =\ast$ lies in $B(N-1)$ so $\gamma=sx$ lies in $B(N)$ and
$\gamma=s$ gives a contradiction.

If $\gamma=xs$ then $w_2x$ evaluates to a subgroup element.  If $x$ is
a stable letter other than $s^{-1}$ then we have the previous case.
So $w_2x\in Z^*$, and there are two further subcases to consider,
depending on $|w_2|$.

\subsubsection*{Case 2b}
If $|w_2| \leq \frac{k}{2}$ then the path that retraces $w$ to
$\overline{w_1}$ then travels along the bottom of the strip to $g'$,
shown in bold in Figure \ref{figA}(b), lies inside $B(N)$ and has
length at most $k+2$.

\subsubsection*{Case 2c}
If $|w_2|>\frac{k}{2}$ then we can argue as follows.  Let $v$ be the
geodesic word in the associated subgroup for $w_2x$ and let $z$ be the
last generator $u_i$ (or $v_i$) of $v$, and let $z_2$ be the last
letter in $Z$ of $z$. That is, $v=v'z=v'z_1z_2$.  The path
$w_2xz_2^{-1}$ is not geodesic in $(A,Z)$, so applying the \fftp\ at
most twice we find a path $y$ of length at most $|w_2|$ which
$2k$-fellow travels it.(That is, $w_2xz_2^{-1}$ and $y$ $k$-fellow
travel an intermediate path $q$.  See Figure \ref{figAc}.
\begin{figure}[ht!]
  \begin{center}
            \subfigure{
      \includegraphics[width=5cm]{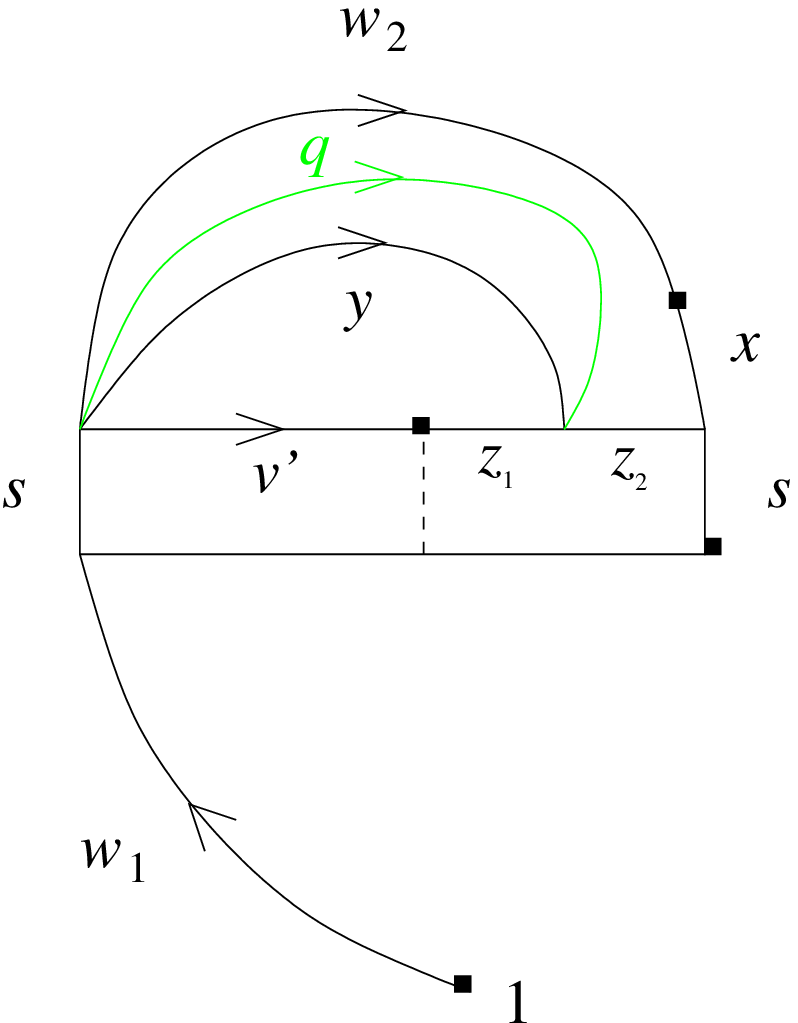}}
	\qquad \qquad
    \subfigure{
      \includegraphics[width=5cm]{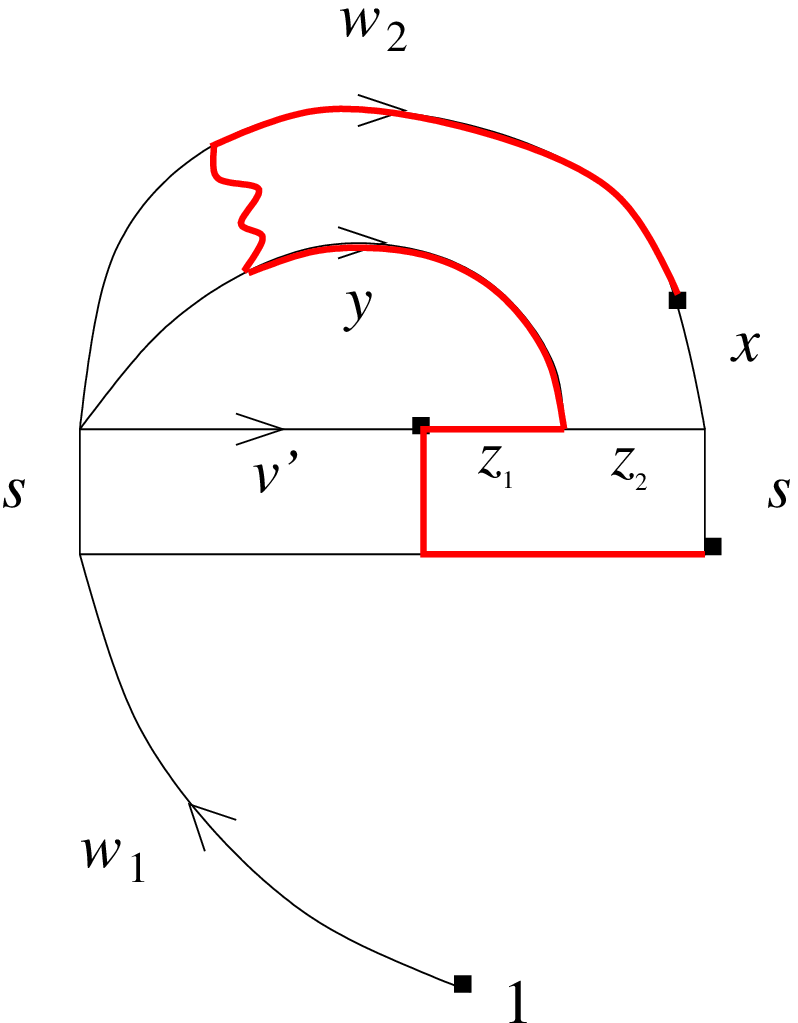}}
   \end{center}
  \caption{Case 2c}
 \label{figAc} 
\end{figure}
Then $y$ lies in $B(N)$ since it starts at the point $\overline{w_1s}
\in B(N-|w_2|)$. We can now find a path connecting $g$ and $g'$ inside
$B(N)$ as follows: Retrace $w$ back to $w(N-\frac{k}{2})$, then across
to $w_1sy(N-\frac{k}{2})$, along $y$ to its end, then around the last
relator of the strip, shown in bold in Figure \ref{figAc}, lies in
$B(N)$ and has length at most $3k + 2\max\{|u_{i_j}|:u_{i_j}
\mathrm{is} \; \mathrm{a} \; \mathrm{generator} \; \mathrm{of}\; U_i
\; \forall i \}$ $\leq 5k$.

\subsection*{Case 3}
$\gamma$ has no stable letters and $w$ has a stable letter.  Then
$w,u$ have parallel stable letter structure.  Let $s\in \{s_i^{\pm
1}\}$ be the last stable letter of $w$, so $w=w_1sw_2, u=u_1su_2$, and
$w_2\gamma u_2^{-1}\in Z^*$ and evaluates to an associated subgroup
element.

\subsubsection*{Case 3a}
If $|w_2|, |u_2| \leq \frac{k}{2}$ then the path that retraces $w$ to
the last strip, then runs across the bottom of the strip, then back
along $su_2$, shown in bold in Figure \ref{figM1}, lies within $B(N)$
and has length at most $4k+2$.  \begin{figure}[ht!]  \begin{center}
\includegraphics[width=10cm]{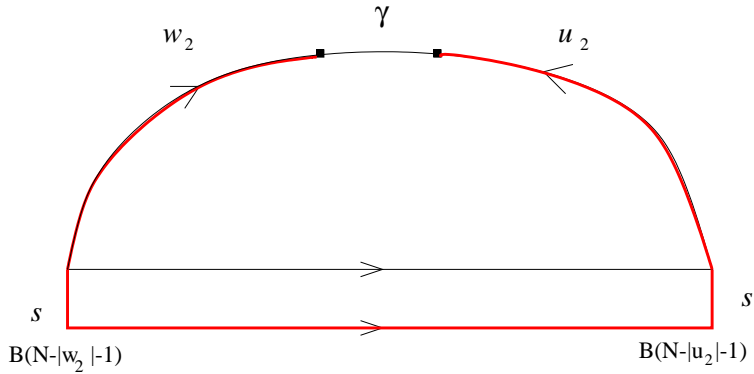} \end{center} \caption{Case 3a:
$|w_2|,|u_2| \leq \frac{k}{2}$.}  \label{figM1}
\end{figure}

Otherwise at least one of $|w_2| > \frac{k}{2}$ or $|u_2| >
\frac{k}{2}$. Without loss of generality assume $|w_2| \geq |u_2| $,
so $|w_2| > \frac{k}{2}$.

\subsubsection*{Case 3b}
If $w_2\gamma u_2^{-1}$ is geodesic then it may not necessarily be a
word in the associated subgroup, that is, it may not run along the top
of the last strip. Let $v$ be a geodesic for it in the associated
subgroup. Let $z$ be the last $|\gamma |$ letters of $v$, so $v=v'z$,
and $|v'| = |w_2|+|u_2|$.  The word $w_2\gamma u_2^{-1}z^{-1}$ has
length $|w_2|+|u_2|+ |\gamma|+|z|$ so is not geodesic, so applying the
\fftp\ in $(A,Z)$ at most four times we can find a geodesic path $y$
for $w_2\gamma u_2^{-1}z^{-1}$ which $4k$-fellow travels it (with
three intermediate paths).

Now since patterns are well behaved, the point $\overline {w_1sv'}$
labeled by $\ast$ in Figure \ref{figM2} must lie in $B(N-|u_2|)$, thus
the path $y$ lies in $B(N)$.
\begin{figure}[ht!]
  \begin{center}
      \includegraphics[width=12cm]{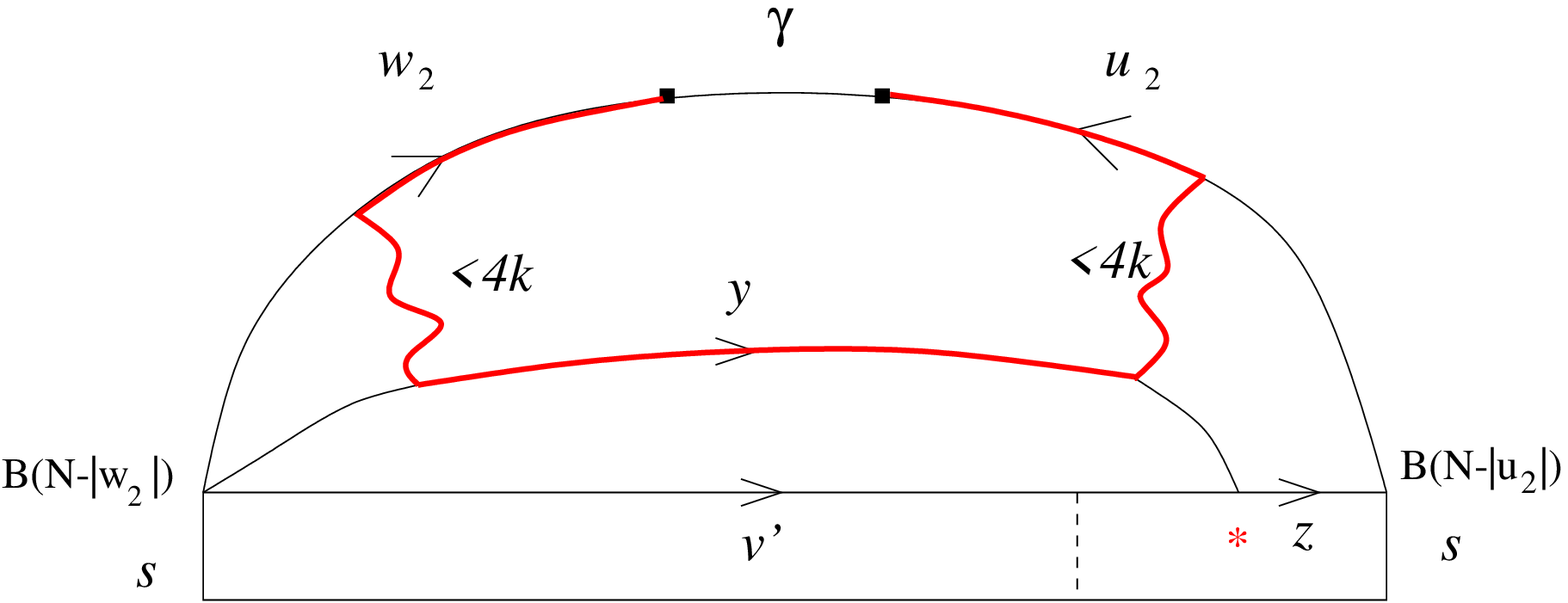}
  \end{center}
  \caption{Case 3b: $|w_2| > \frac{k}{2}, w_2\gamma u_2$  geodesic.}
  \label{figM2}
\end{figure}
The path that retraces $w$ to $w(N-\frac{k}{2})$, then across to
$y(N-\frac{k}{2})$, then along $y$ to $y(N+\frac{k}{2}+|\gamma |)$,
then across to $u(N-\frac{k}{2})$ and back to $g'$ along $u$, shown in
bold in Figure \ref{figM2}, lies within $B(N)$ and has length at most
$10k+2$.  If $|u_2|<\frac{k}{2}$ then the path would follow $y$ until
its end, go around the last relator of the strip (which contains $z$)
and along $u_2$.

\subsubsection*{Case 3c}
If $w_2\gamma u_2^{-1}$ is not geodesic then by the \fftp\ in $(A,Z)$
there is a path $y$ which $k$-fellow travels it.  If $y$ is geodesic,
we can repeat the previous argument to obtain a path of length at most
$10k+2$ inside $B(N)$.  The details are omitted.

\subsubsection*{Case 3d}
If $y$ is also not geodesic then applying the \fftp\ once more we get
a path $p$ of length at most $|w_2|+|u_2|+|\gamma|-2 \leq
|w_2|+|u_2|$, which runs from $\overline {w_1s}$ to $\overline {u_1s}$
so lies in $B(N)$.  Then the path that retraces $w$ to
$w(N-\frac{k}{2})$, then across to $y(N-\frac{k}{2})$, then across to
$p(N-\frac{k}{2})$, along $p$ to $p(N+\frac{k}{2}+|\gamma|)$, over to
$u(N-\frac{k}{2})$ and back to $g'$, shown in bold in Figure
\ref{figM3}, lies within $B(N)$ and has length at most $6k+2$.  If
$|u_2|<\frac{k}{2}$ then the path would follow $p$ until its end, then
run along $u_2$.
\begin{figure}[ht!]
  \begin{center}
      \includegraphics[width=9.7cm]{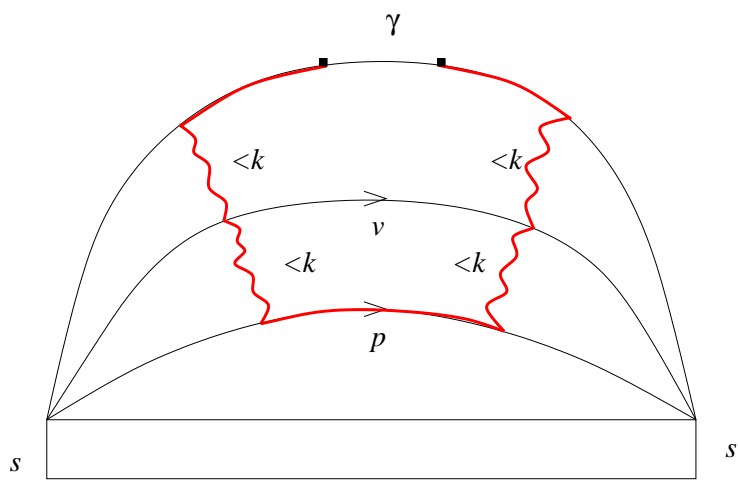}
  \end{center}
  \caption{Case 3d: $|w_2| > \frac{k}{2}, w_2\gamma u_2, y$ not geodesic.}
  \label{figM3}
\end{figure}

In each case the maximum length of a path needed to connect $g$ to
$g'$ inside $B(N)$ is $10k+2$, and the result follows.
\end{proof}

\bibliography{refs}

\begin{thebibliography}{10}

\bibitem{BB}
Noel Brady and Martin Bridson.
\newblock On the absense of biautomaticity for graphs of abelian groups.
\newblock Unpublished.

\bibitem{MR1744486}
Martin~R. Bridson and Andr{\'e} Haefliger.
\newblock {\em Metric spaces of non-positive curvature}.
\newblock Springer-Verlag, Berlin, 1999.

\bibitem{MR88a:20049}
James~W. Cannon.
\newblock Almost convex groups.
\newblock {\em Geom. Dedicata}, 22(2):197--210, 1987.

\bibitem{Ethesis}
Murray Elder.
\newblock Automaticity, almost convexity and falsification by fellow traveler
  properties of some finitely generated groups.
\newblock PhD Dissertation, University of Melbourne, 2000.

\bibitem{MR2022477}
Murray~J. Elder.
\newblock A non-{H}opfian almost convex group.
\newblock {\em J. Algebra}, 271(1):11--21, 2004.

\bibitem{MR93i:20036}
David B.~A. Epstein, James~W. Cannon, Derek~F. Holt, Silvio V.~F. Levy,
  Michael~S. Paterson, and William~P. Thurston.
\newblock {\em Word processing in groups}.
\newblock Jones and Bartlett Publishers, Boston, MA, 1992.

\bibitem{MR94j:20043}
S.~M. Gersten.
\newblock The automorphism group of a free group is not a cat(0) group.
\newblock {\em Proc. Amer. Math. Soc.}, 121(4):999--1002, 1994.

\bibitem{MR83j:68002}
John~E. Hopcroft and Jeffrey~D. Ullman.
\newblock {\em Introduction to automata theory, languages, and computation}.
\newblock Addison-Wesley Publishing Co., Reading, Mass., 1979.
\newblock Addison-Wesley Series in Computer Science.

\bibitem{MR58:28182}
Roger~C. Lyndon and Paul~E. Schupp.
\newblock {\em Combinatorial group theory}.
\newblock Springer-Verlag, Berlin, 1977.
\newblock Ergebnisse der Mathematik und ihrer Grenzgebiete, Band 89.

\bibitem{MR96c:20066}
Walter~D. Neumann and Michael Shapiro.
\newblock Automatic structures, rational growth, and geometrically finite
  hyperbolic groups.
\newblock {\em Invent. Math.}, 120(2):259--287, 1995.

\bibitem{MR96m:20058}
Daniel~T. Wise.
\newblock A non-{H}opfian automatic group.
\newblock {\em J. Algebra}, 180(3):845--847, 1996.

\end{thebibliography}
\bibliographystyle{plain}

\end{document}